\documentclass[11pt, a4paper, twocolumn]{article} % 10pt font size (11 and 12 also possible), A4 paper (letterpaper for US letter) and two column layout (remove for one column)
\usepackage{graphicx}
\usepackage{subcaption}
\usepackage[english]{babel} % English language hyphenation

\usepackage{microtype} % Better typography

\usepackage{amsmath,amsfonts,amsthm} % Math packages for equations

\usepackage[svgnames]{xcolor} % Enabling colors by their 'svgnames'

\usepackage{booktabs} % Horizontal rules in tables

\usepackage{lastpage} % Used to determine the number of pages in the document (for "Page X of Total")

\usepackage{graphicx} % Required for adding images

\usepackage{enumitem} % Required for customising lists
\setlist{noitemsep} % Remove spacing between bullet/numbered list elements

\usepackage{sectsty} % Enables custom section titles
\allsectionsfont{\usefont{OT1}{phv}{b}{n}} % Change the font of all section commands (Helvetica)

\usepackage{geometry} % Required for adjusting page dimensions

\usepackage[T1]{fontenc} % Output font encoding for international characters
\usepackage[utf8]{inputenc} % Required for inputting international characters

\usepackage{XCharter} % Use the XCharter font

\usepackage{fancyhdr} % Needed to define custom headers/footers
\fancypagestyle{firstpage}{ % Page style for the first page with the title
	\fancyhf{}
}

\newcommand{\institution}[1]{{\footnotesize\usefont{OT1}{phv}{m}{sl}\color{Black}#1}} % Institutions style (Helvetica)

\usepackage{titling} % Allows custom title configuration

\newcommand{\HorRule}{\color{DarkGoldenrod}\rule{\linewidth}{1pt}} % Defines the gold horizontal rule around the title

\pretitle{
	\vspace{-30pt} % Move the entire title section up
	\HorRule\vspace{10pt} % Horizontal rule before the title
	\fontsize{32}{36}\usefont{OT1}{phv}{b}{n}\selectfont % Helvetica
	\color{DarkRed} % Text colour for the title and author(s)
}

\posttitle{\par\vskip 15pt} % Whitespace under the title

\preauthor{} % Anything that will appear before \author is printed

\postauthor{ % Anything that will appear after \author is printed
	\vspace{10pt} % Space before the rule
	\par\HorRule % Horizontal rule after the title
	\vspace{20pt} % Space after the title section
}

\usepackage{lettrine} % Package to accentuate the first letter of the text (lettrine)
\usepackage{fix-cm}	% Fixes the height of the lettrine

\usepackage{xstring} % Required for string manipulation

\newcommand{\AuthorNames}{Eleuterio F.  Toro}
\newcommand{\ShortTitle}{Personal Memories of Philip L.  Roe}

\title{\LARGE   Philip L.  Roe,  the British Scientist.\\ 
The Cranfield Years and Beyond} % The article title

\author{\sc Eleuterio F.  Toro
	\newline\newline % Space before institutions
	 \institution{\sc Professor Emeritus,  Mathematics, University of Trento, Italy\\
	 Life Fellow of Clare Hall, University of Cambridge, UK\\
	 \tt Personal webpage: https://eleuteriotoro.com \\
	  Email: eleuterio.toro@unitn.it}}
	 
\date{} % Add a date here if you would like one to appear underneath the title block, use \today for the current date, leave empty for no date

\usepackage{fancyhdr}
\begin{document}

\maketitle

\thispagestyle{firstpage}

%\twocolumn[%
  %\begin{center}
    %\begin{minipage}{0.95\textwidth}
      %\small % rende il blocco più compatto
     % \selectlanguage{italian}
      %\textbf{Sommario}\par
      %Si presenta la Scuola del GNFM, che quest’anno celebra il suo cinquantenario. 
      %I temi affrontati e i prestigiosi relatori riflettono %l’evoluzione delle ricerche 
      %e il continuo rinnovamento della Fisica Matematica in %Italia.
      %\vskip0.8\baselineskip % spazio contenuto tra sommario %e abstract

      %selectlanguage{english}
      %\textbf{Abstract}\par
      %This paper presents the GNFM School, which this year %celebrates its fiftieth anniversary. 
      %The topics covered and the distinguished speakers %reflect the evolution of research 
      %and the ongoing development of Mathematical Physics in Italy.
    %\end{minipage}
  %\end{center}
  %\vskip0.5\baselineskip % piccolo spazio dopo il blocco
%]

%----------------------------------------------------------------------------------------

%PENDING: Alan's feedback; final read; submit.

\begin{center}
\begin{minipage}{0.95\columnwidth}
%\small
\noindent\textbf{\small{ABSTRACT}} 

{{
Philip L. Roe, the renowned British scientist, passed away on April 26, 2026, at the age of 87, while Professor Emeritus at the University of Michigan.

My association with Phil Roe began in 1984, when he joined the Cranfield Institute of Technology, UK, as Professor of Computational Fluid Dynamics, while I was a postdoctoral research fellow at the same institution.
For more than four decades, his incisive ideas,  physical principles and distinctive approach to research on numerical methods for solving hyperbolic partial differential equations have had a lasting influence on my own approach to research in the field.

Here I recount some personal memories and anecdotes emerging from the rewarding experience of interacting with Phil during the Cranfield years and beyond.  I also review—quite succinctly in places—some of the main research topics that he pursued throughout much of his academic career.

These notes are intended as a modest homage to Roe's impact on the subject, as well as to his influence on the careers of generations of scientists in both academia and industry.
}}
\end{minipage}
\end{center}

%----------------------------------------------------------------------------------------
%	ARTICLE CONTENTS
%----------------------------------------------------------------------------------------

\section{Roe's Cranfield Years}

%DONE

 I begin by recalling some historical notes, recounting personal anecdotes, and briefly describing Roe’s contributions to research on numerical methods and postgraduate teaching at Cranfield.

\subsection{From RAE to Cranfield}

%DONE

The year 1984 marked a turning point in the professional life of Phil Roe. He moved from the Royal Aircraft Establishment (RAE) in Bedford to the Cranfield Institute of Technology—referred to hereafter as Cranfield—located near Cranfield Village. Although his new workplace was farther away—requiring a 16 km drive west instead of his usual 10 km journey north—he chose to keep his family home in Bedford.

The adjustment was mostly straightforward. However, on more than one occasion—he once told me—he found himself arriving at his former workplace out of habit.

What are you doing here, Phil?
his colleagues would ask, as he turned around and headed back to his car.

RAE at Bedford was a UK government laboratory dedicated to advanced aeronautical research. Work at RAE was underpinned by major wind tunnel facilities that enabled sophisticated aerodynamic testing, forming the core of critical investigations into high-speed flight, stability, and control. At that time, these capabilities provided the experimental foundation for landmark aircraft programmes such as Concorde and the Harrier, supporting the UK’s leadership in both supersonic and vertical/short take-off and landing technologies.

While at RAE, Phil Roe focused primarily on experimental aerodynamics, contributing to the understanding of complex flow behaviour through wind tunnel research. But in his final eight years at RAE, Roe's work shifted towards the development of computational algorithms for the theoretical simulation of fluid dynamics. This emerging field—later known as Computational Fluid Dynamics (CFD)—became a cornerstone of modern aerospace engineering and beyond, complementing experimental methods and enabling theoretical prediction of aerodynamic performance.

CFD became Roe's passport to academic research and postgraduate teaching at Cranfield. This UK higher education institution started as the College of Aeronautics in 1946, founded on the former Royal Air Force base, RAF Cranfield. Then, in 1969, it was granted university status by Royal Charter, with the addition of more departments such as mechanical engineering, materials, applied mathematics and management; in 1993, it was renamed Cranfield University.

Cranfield, at the Cranfield campus, was exclusively a postgraduate university, offering a wide range of master's and PhD degrees. Cranfield's undergraduate campuses were located at Shrivenham, west of Oxford, and Silsoe, south of Bedford.  

The move from RAE to Cranfield may at first sight appear to mark a drastic technical change in Roe's professional life. As a matter of fact, however, his experimental aerodynamics experience was entirely functional to the theoretical and computational modelling of compressible fluid dynamics. His approach to CFD was strongly, and positively, influenced by experimental research; this was especially reflected in his teaching of postgraduate aeronautical engineers. This was, and still is, a point of the utmost importance in keeping theoretical research and computational simulations close to physical reality. This is a principle that applies to all sciences and is particularly relevant in educating young scientists in an era of an increasing number of new theoretical research tools. Phil Roe was a passionate advocate of this principle.

I joined Cranfield in 1983, one year before Roe, as a postdoctoral research fellow. My research background was in finite element methods for elliptic free-boundary problems, while my new remit was mathematical and computational modelling of propulsion systems under the guidance of Professor John F.  Clarke (1927–2013).

\subsection{Getting to know my driver}

%DONE

I confess I did not, initially, fully appreciate the significance of Roe's arrival at Cranfield. But that would soon change. Shortly after joining our department, Phil kindly invited me to accompany him to a one-day workshop at the University of Cambridge, just over an hour’s drive from Cranfield. The meeting focused on numerical methods for hyperbolic equations, a topic that was almost entirely new to me.

We set off early that morning in Phil’s large, old Volvo, making our way slowly across the countryside. We chatted about a variety of general topics, as we did not yet have much in common in technical terms. He was quite familiar with Chile, my country of origin, and told me about the pioneering efforts of President Salvador Allende in the early 1970s to introduce high technology into business and the national economy. This was the Cybersyn Project, led by the British operations research scientist Stafford Beer (1926–2002) \cite{Beer:1995a} and the Chilean engineer Fernando Flores. I was impressed by Phil’s breadth of knowledge about my country.

The workshop itself was relatively small. The participants seemed friendly, but I must confess that I remember almost nothing about the presentations. The most interesting part of the gathering took place during the coffee break, as is often the case at scientific meetings of this kind.

I struck up a conversation with a young man named Reinald L\"{o}hner. We quickly discovered that we had something in common: the Spanish language. L\"{o}hner was Argentinian  and had just completed his PhD at Swansea University under the supervision of Professor O. Zienkiewicz (1921–2009). Fortunately, I was already aware of Zienkiewicz’s scientific standing. Interestingly, Reinald’s external PhD examiner had been Phil Roe,  he told me. Perhaps noticing my inquisitive expression, he explained to me, in clear Spanish, that Phil Roe was a highly renowned figure in the field of numerical methods for solving partial differential equations (PDEs) of hyperbolic type.  Reinald L\"{o}hner was later to become a very prominent academic.

On the way back to Cranfield that evening, driving through the dark, I found myself casting furtive glances at the man at the wheel to my right. I felt fortunate to be driven by such a distinguished chauffeur.

\subsection{Phil's office at Cranfield}

%DONE

Phil had—quite rightly—an enormous office, at least by Cranfield standards. Rumour had it that some postdocs, housed in {\it temporary} portacabins, thought it excessive. Yet, on entering the room, what struck you was not its size but its centre of gravity: Phil himself, seated in a comfortable chair in one corner. Behind a large radiator under the windowsill lay the noisy airfield; before him, an imposing yet curiously {\it invisible} desk, buried beneath a vast accumulation of journal papers and technical reports. That corner was Roe's true living space.

\begin{figure}[h]
      \centerline{
             \includegraphics[scale=0.14, angle=0]{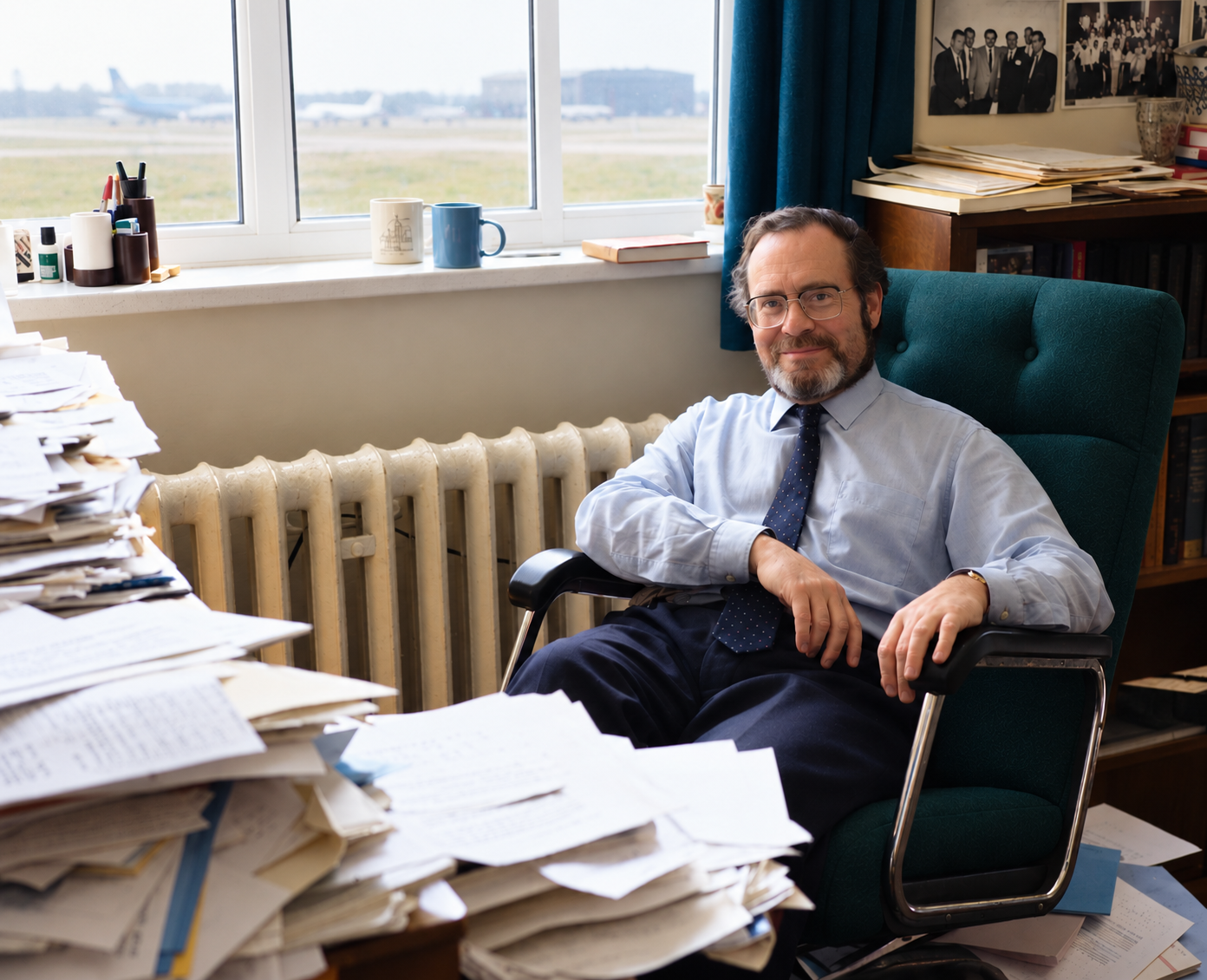}
      }
      \caption{Professor Philip L.  Roe in his office at the Cranfield Institute of Technology, UK (1987 circa).}
      \label{fig:Phil-living-space}
\end{figure}

The rest of the room—absent from  Fig. \ref{fig:Phil-living-space}—was a large, rectangular expanse dominated by an equally large meeting table, placed with near-geometric symmetry at its centre and surrounded by eight or ten chairs. The table’s surface, several of the chairs, and even parts of the surrounding floor were occupied by still-sealed boxes, presumably filled with books and papers. Navigating this uninhabited territory was no small feat.

I cannot recall whether those boxes were ever opened or even moved after Phil’s arrival at Cranfield from RAE in 1984. It is entirely possible that they remained untouched until he left for Michigan in 1990. In any case, the modest corner he inhabited seemed more than sufficient for his needs—an arrangement that reflected the distinctive traits evident in his daily conduct: humility and minimalism. Take only what you need, and no more.

Within that corner, the non-static vast accumulation of papers atop his desk would, most of the time, exhibit a fairly steady arrangement of three contiguous piles. These formed a configuration, just in front of Phil, that would have delighted any numerical analyst. They defined, with remarkable clarity, the initial conditions for two Riemann problems—those arising at the left and right interfaces of a finite volume in space-time. The piles corresponded to three constant cell averages for three consecutive cells, labelled $i-1$ (left), $i$ (centre), and $i+1$ (right). The middle pile—cell $i$—was the lowest, just above the level of the desk: the vacuum state. It was there that Phil carried out his characteristically tidy algebra, working steadily in a large copybook.

Yet these piles were far more than abstract constructs. They formed part of a treasure trove of knowledge—housing papers not easily found in any library. I would often visit him in search of a particular article I had failed to locate elsewhere. I recall that on one occasion, after I mentioned the author, title, and year, Phil leaned back, momentarily pensive. His gaze settled on a seemingly arbitrary point amid the surrounding mountains of paper. It appeared an impossible task—to retrieve a specific document from such apparent chaos within any reasonable time.

Then, without hesitation, he leaned forward, extended both arms, and extracted the desired paper as if performing a conjuring trick. What had seemed disorder to me was, to Phil, in fact, a highly personal and perfectly functional filing system.

\subsection{Cranfield's airfield}

%DONE

During my time at Cranfield,  the College of Aeronautics was  one of Cranfield's main academic structures.
It comprised three departments: Structures, Aviation, and Aerodynamics.  Phil Roe, John Clarke, and I were all members of the Aerodynamics Department. The vast majority of the academic staff were engineers—predominantly aeronautical engineers—for whom aircraft were not just a profession but a lifelong passion. There were a few physicists and some mathematicians, for whom the excitement of aircraft operations, noisy helicopters, wind tunnels, and shock tube laboratories was perhaps less compelling.

One day, John Stollery, an experimental aerodynamicist and Head of the Aerodynamics Department, called eight of us to a meeting in his office. Phil Roe,  an aeronautical engineer by trade,  was present, as was my colleague John Clarke, trained as an aeronautical engineer but working in practice as an applied mathematician. I felt somewhat out of place—both as the most junior member and as a mathematician in such a strongly aeronautical environment.  Stollery’s office occupied a corner position on the second floor, offering a wide view over the busy airfield through large windows—clearly a vantage point he appreciated.

On what was a rare sunny afternoon at Cranfield, Stollery drew the curtains to shield us from the glare as the meeting began. Despite the closed windows and heavy curtains, the noise from the airfield—especially the helicopters—was intense, and I found it difficult to concentrate. My colleagues, however, seemed entirely absorbed.

At one point, the discussion abruptly halted. Everyone fell silent, then turned their heads as if trying to identify a particular sound amid the general noise. Someone exclaimed,

—It’s the XYZ123!—

Instantly, they rushed to the windows, pulling back the curtains and scanning the sky. A few even stood on chairs to get a better view. I remained seated—the only one unmoved by the excitement. Clearly, I was not an aeronautical engineer. I cannot recall exactly where Phil Roe stood at that moment, but I suspect he maintained rather more composure than most.

I clarify that I do not remember the exact designation for that special aircraft. 

\begin{figure}[h]
      \centerline{
             \includegraphics[scale=0.14, angle=0]{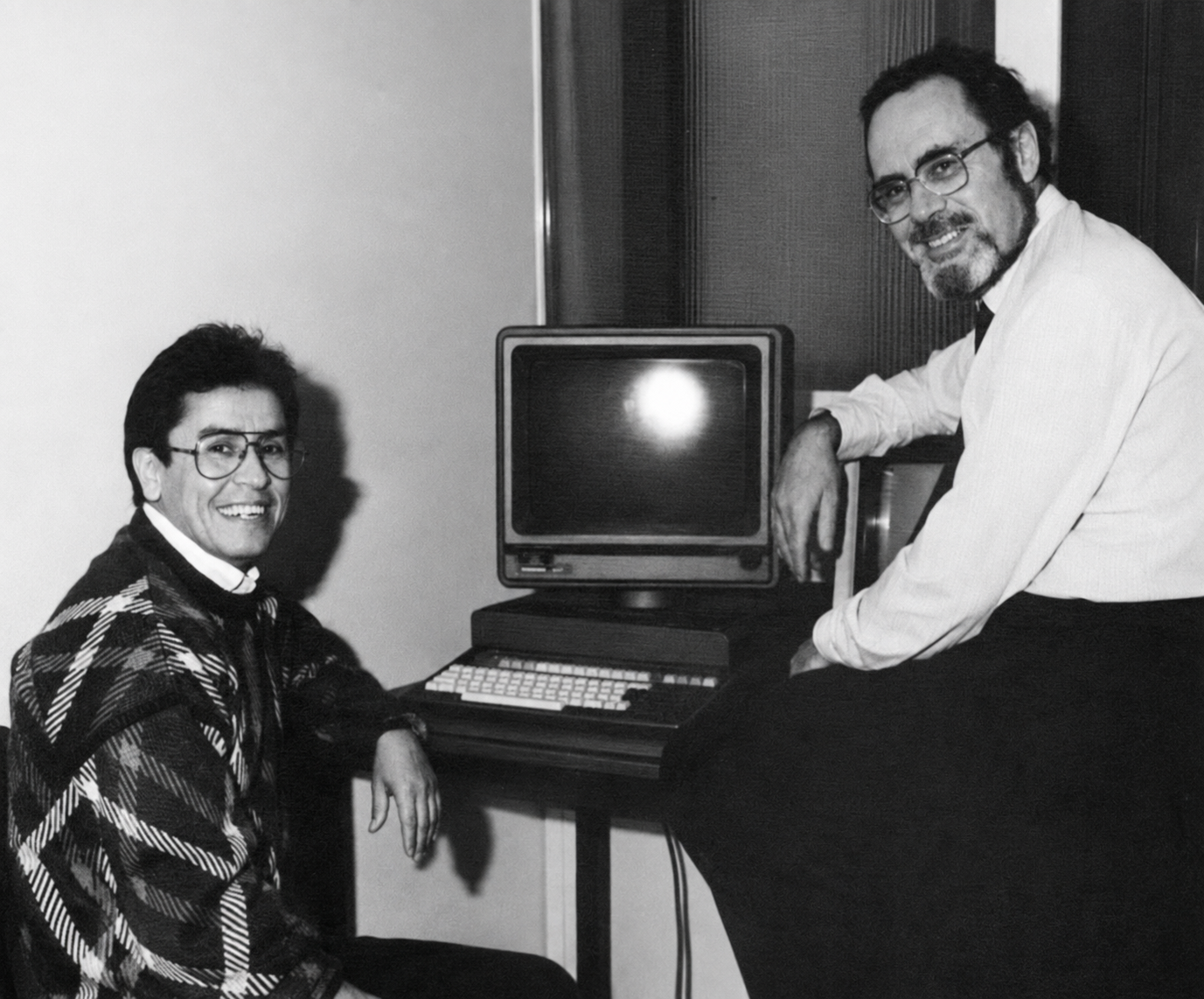}
      }
      \caption{Professor Philip  L. Roe (right) and Dr Eleuterio F. Toro.   College of Aeronautics,  Cranfield, UK (1987 circa).}
      \label{fig:PLRoe-EFToro}
\end{figure}

\subsection{The Cranfield CFD group}

%DONE

The arrival of Phil Roe at Cranfield in 1984 marked a turning point in postgraduate teaching and research within the Aerodynamics Department. At the time, there was already some activity in the mathematical modelling of fluid dynamics processes—particularly combustion—under the leadership of John Clarke, an expert in asymptotic analysis. A notable innovation introduced by Phil was a new module in the MSc Aerodynamics programme entitled Computational Fluid Dynamics. I attended most of his lectures and found them highly instructive.

Cranfield had a well-established tradition of offering short MSc/PhD-level courses to external participants. These courses were professionally organised and contributed to the Department’s external income. A dedicated administrator managed their commercial delivery.

Under Phil Roe’s leadership, a short course in CFD was introduced. John Clarke contributed lectures on gas dynamics and combustion, and in time I also became involved, delivering lectures on exact Riemann solvers for the Euler equations and Glimm’s Random Choice Method. Participants came from both academia and industry, primarily from the UK and continental Europe, but also from other continents. The CFD course proved to be a great success and had a lasting influence on my own academic work. I remember that, for promotional purposes, Phil devised the playful rebranding of Cranfield as
{\bf C}ran{\bf F}iel{\bf D},
highlighting the acronym CFD.

Research activity in CFD expanded with the arrival of new PhD students,  including Smadar Karni,  James Quirk,  Lloyds G.  Simmonds and H. C.  Lin, all supervised by Phil Roe.  John Clarke supervised, with some contributions from me,  PhD students Alan Dawes,  Gurlat Singh and Nikolaos Nikiforakis.  Another important addition to the group was Jack Pike,  who joined Cranfield as a research fellow.  Pike had previously worked with Roe at RAE on Riemann solvers, and while at Cranfield he communicated new contributions on the subject \cite{Pike:1993a}.  All these developments underwent a significant shift when Phil left for Michigan. Nevertheless, strong foundations had been laid.  

After Phil left for Michigan in 1990,  I was promoted to a lectureship in CFD, after several years in a permanent position as Senior Research Fellow. The invaluable experience gained during Phil’s time at Cranfield enabled me to establish a new MSc programme in CFD in 1990/1991. This programme was highly successful and, in various forms, has continued at Cranfield to the present day, in 2026. At the time, it attracted outstanding graduates, predominantly in engineering, physics, and mathematics. Notable MSc students included D. Elsworth, A. Chakraborty, I. R. Dawkins, C. J. Griesel, R. Brown, P. S. Darling, J. Orozco, Andrew Smith, D. N. Cheney, P. Burton, N. J. Hughes, I. D. Fraser, D. Brocklehurst, C. C. Chou, Y. Mee-Sin, and M. Cáceres. PhD students included N. Nikiforakis and C. Lowe (both co-supervised with J. F. Clarke), while S. J. Billett and E. P. Boden were entirely under my supervision.

In collaboration with W. S. Speares and M. Spruce, we developed the HLLC Riemann solver for the time-dependent Euler equations 
\cite{Toro:1994c}.  In parallel, with A. Chakraborty, we formulated an HLLC scheme for the 2D steady supersonic Euler equations \cite{Toro:1994a}. N. Nikiforakis conducted a theoretical study of detonation waves using the then-recent high-resolution WAF method \cite{Toro:1989a}. Together with S. J. Billett, we reformulated the WAF approach in a multidimensional setting, building on its original formulation for one-dimensional conservation laws \cite{Billett:1994a}, \cite{Billett:1997a}. With MSc students M. C\'aceres and D. N. Cheney, we also implemented and tested very early ideas on the ADER high-order methodology in the early 1990s.

These developments, even if indirectly, can be traced back to the research foundations laid during Roe’s time at Cranfield.

Phil Roe made remarkable contributions to the field of Computational Fluid Dynamics.  A personal account by Roe himself is found in \cite{Roe:2020a}.  In what follows I briefly review some of his early contributions,  with which I am more familiar. The review is not exhaustive. 

\section{Roe's Riemann Solver}

If any single work best illustrates Roe’s impact on CFD,  it is his 1981 paper {\it Approximate Riemann solvers, parameter vectors, and difference schemes} \cite{Roe:1981a}.  

A web search confirms it to be the single most-cited article in the field, ahead of Godunov’s seminal paper \cite{Godunov:1959a}; at the time of writing, Google Scholar listed 14876 citations,  out of a total of 32263 citations (April 2026); that is,  nearly half of all Roe's citations in Google Scholar are on his single 1981 paper.

\subsection{The finite volume framework}
\label{sec:FV}

\begin{figure}
         \centerline{
           \includegraphics[scale=0.7]{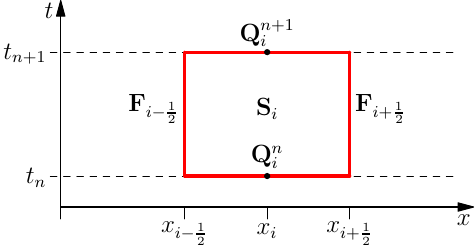}  
           }
           \caption{Representation of finite volume formula (\ref{eq:Roe2}) in $x$-$t$ space.}
           \label{fig:FVformula}
\end{figure}

Roe's method is best understood in the framework of finite volume methods. 
Consider a general system of $m$  hyperbolic balance laws
\begin{equation}                                                  \label{eq:Roe1}
      \partial_{t}{\bf Q}(x,t) +  \partial_{x}{\bf F}({\bf Q}(x,t))={\bf S}({\bf Q}(x,t)) \;
\end{equation}
and the control volume   $V_{i}^{n}=[x_{i-\frac{1}{2}}, x_{i+\frac{1}{2}}] \times [t_{n},t_{n+1}]$ in $x$-$t$ space, as illustrated in Fig.   \ref{fig:FVformula}.
Integration of (\ref{eq:Roe1}) in $V_{i}^{n}$ gives
\begin{equation}                                                     \label{eq:Roe2}
    \hat{\bf Q}_{i}^{n+1}=\hat{\bf Q}_{i}^{n} - \frac{\Delta t}{\Delta x}(\hat{\bf F}_{i+\frac{1}{2}}-\hat{\bf F}_{i-\frac{1}{2}} ) + \Delta t \hat{\bf S}_{i} \;.
\end{equation}
Formula (\ref{eq:Roe2}) is {\bf  exact}, with definitions
\begin{equation}                                                  \label{eq:Roe3}
      \left.\begin{array}{l}
      \displaystyle{ \hat{\bf Q}_{i}^{n}=\frac{1}{\Delta x} \int_{x_{i-\frac{1}{2}}}^{x_{i+\frac{1}{2}}} {\bf Q}(x,t_{n}) dx} \;, \\[13pt]
      \displaystyle{\hat{\bf F}_{i+\frac{1}{2}} = \frac{1}{\Delta t} \int_{t_{n}}^{t_{n+1}} {\bf F}( {\bf Q}(x_{i+\frac{1}{2}}, t) ) d t} \;, \\[10pt]
      \displaystyle{\hat{\bf S}_{i} = \frac{1}{\Delta t \Delta x} \int_{t_{n}}^{t_{n+1}} \int_{x_{i-\frac{1}{2}}}^{x_{i+\frac{1}{2}}} {\bf S}( {\bf Q}(x,t) ) d x d t} \;.
     \end{array} \right\} 
\end{equation}
Formula (\ref{eq:Roe2}) may be interpreted as a one-step numerical scheme,  if the  space and time integrals in (\ref{eq:Roe3}) are appropriately computed. In this case formula (\ref{eq:Roe2}) becomes a numerical method for evolving spatial cell integral averages $\hat{\bf Q}_{i}^{n}$,  with numerical flux $\hat{\bf F}_{i+\frac{1}{2}}$ and numerical source $\hat{\bf S}_{i}$.  Fig.  \ref{fig:FVformula} illustrates scheme (\ref{eq:Roe2}). 
For the rest of this section we assume the homogeneous case in (\ref{eq:Roe1}),  that is ${\bf S}({\bf Q}(x,t))={\bf 0}$.  

Back in 1959 Godunov \cite{Godunov:1959a} suggested to define the numerical  flux ${\bf F}_{i+\frac{1}{2}}$ at the cell interface by using the solution of the Riemann problem
\begin{equation}     \label{eq:Roe4}
 \left.\begin{array}{l}
 \partial_t {\bf Q} + \partial_x {\bf F}({\bf Q}) = {\bf 0} \;, \hspace{2mm}  x \in (-\infty, \infty) \;,\hspace{2mm} t>0 \;, \\
  {\bf Q}(x,0) = \left\{
       \begin{array}{ccc}
           {\bf Q}_{i}^{n} = {\bf Q}_{L}     & \mbox{   if   } & x <  0 \;,\\
           {\bf Q}_{i+1}^{n}={\bf Q}_{R} & \mbox{   if   } & x >  0 \;.
       \end{array}\right.
\end{array}\right\}   
\end{equation}
In the original paper of 1959,  Godunov solved the Riemann problem (\ref{eq:Roe4}) exactly \cite{Godunov:1959a},  to get the similarity solution $\tilde{\bf Q}(x/t)$,  in local coordinates $x$, $t$.  After appropriate sampling the constant value along the interface $x/t=0$ is found and from (\ref{eq:Roe3}) the numerical flux becomes
\begin{equation}    \label{eq:Roe5}
    {\bf F}_{i+\frac{1}{2}}= {\bf F}(\tilde{\bf Q}(0))\;.
\end{equation}

\subsection{Roe's linearisation}

 A key point of Roe's approach to designing a numerical flux for Eq.  (\ref{eq:Roe2})
was based on the  argument  {\it that in existing schemes much of this information is degraded,  and that only certain features of the exact solution are worth striving for}.  

Roe assumed that the nonlinear system in (\ref{eq:Roe4}) could be replaced by a linear system  with constant coefficients leading to the linear Riemann problem:
\begin{equation}     \label{eq:Roe6}
 \left.\begin{array}{l}
 \partial_{t}{\bf Q}+    \tilde{\bf A}({\bf Q}_{L},{\bf Q}_{R}) \partial_{x}{\bf Q}= {\bf 0}\;, \\
  {\bf Q}(x,0) = \left\{
       \begin{array}{ccc}
           {\bf Q}_{i}^{n} = {\bf Q}_{L}     & \mbox{   if   } & x <  0 \;,\\
           {\bf Q}_{i+1}^{n}={\bf Q}_{R} & \mbox{   if   } & x >  0 \;.
       \end{array}\right.
\end{array}\right\}   
\end{equation}
The sought matrix $\tilde{\bf A}$,  often called the {\bf Roe matrix}, is related to the Jacobian matrix  ${\bf A}=\partial {\bf F}/\partial {\bf Q}$ and is expected to satisfy the following  properties:
\begin{enumerate}
\item {\bf Hyperbolicity. } $\tilde{ \bf A }$ is required to have real eigenvalues $\tilde \lambda_i = \tilde \lambda_i ({\bf Q}_{L}, {\bf Q}_{R})$,  ordered as
\begin{equation}                                              \label{eq:Roe7}
  \tilde{\lambda}_1 \leq \tilde{\lambda}_2 \leq \cdots \leq \tilde{\lambda}_m \;,
\end{equation}
and a corresponding complete set of linearly independent right eigenvectors
\begin{equation}                                              \label{eq:Roe8}
  \tilde{\bf K }^{(1)} \;, \hspace{3mm} \tilde{\bf K }^{(2)} \;, \cdots \;, 
  \hspace{3mm} \tilde{\bf  K }^{(m)} \;. 
\end{equation}
\item {\bf Consistency. } The matrix $\tilde{ \bf A }$ is required to be consistent with the exact Jacobian, that is
\begin{equation}                                              \label{eq:Roe9}
  \tilde{\bf A } (\bf  Q,\bf Q ) = \bf A (\bf Q ) \;.
\end{equation}
\item {\bf Conservation. } Across discontinuities  $\tilde{ \bf A }$ is required to satisfy the jump condition
\begin{equation}                                              \label{eq:Roe10}
  {\bf F}({\bf Q}_{R}) - {\bf F}({\bf Q}_{L})
  = \tilde{\bf A } \; ({\bf Q}_{R} - {\bf Q}_{L}) \;. 
\end{equation}
\end{enumerate}

The construction of  matrices satisfying the above properties for {\bf general hyperbolic systems} can be very complicated.  For the {\bf specific} case of the Euler equations of gas dynamics Roe \cite{Roe:1981a} 
proposed a relatively simple way of constructing a matrix $\tilde{\bf A }$.  Later, Roe and Pike \cite{Roe:1984a} proposed a  simpler approach, where the explicit construction of $\tilde{\bf A}$ is avoided.

A breakthrough in constructing $\tilde{\bf A}$ resulted from Roe's original idea of introducing a {\bf parameter vector}  $\bf W$, such that both the vector of conserved variables $\bf Q$ and the flux vector  $\bf F(\bf Q)$ could be expressed 
in terms of $\bf W$. That is
\begin{equation}                                              \label{eq:Roe11}
  \bf Q = \bf Q(\bf W) \;, ~~ \bf F = \bf F(\bf W) \;.
\end{equation}
Two important steps then follow. First, the changes
\begin{equation}                                              \label{eq:Roe12}
  \Delta {\bf Q} = {\bf Q}_{R} - {\bf Q}_{L} \;,  ~~
  \Delta {\bf F} = {\bf F}({\bf Q}_{R}) - {\bf F}({\bf Q}_{L})
\end{equation}
can be expressed in terms of the change $\Delta {\bf W} = {\bf W}_{R}-{\bf W}_{L}$. Then, averages are obtained in terms  of {\em simple arithmetic means} of ${\bf W}$, which may appear counter-intuitive.  

Next, we illustrate the technique as applied to a simple set of hyperbolic conservation laws.

\subsection{Roe's method for a simple system}

Consider the isothermal equations of gas dynamics
\begin{equation}                                              \label{eq:Roe13}
  \left. \begin{array}{c}
     \partial_t {\bf Q} + \partial_x {\bf F}({\bf Q}) = {\bf 0}  \;,    \\[6pt]
    {\bf Q} \equiv 
    \left[ \begin{array}{c}
      q_1       \\
      q_2
    \end{array} \right]
     \equiv 
    \left[ \begin{array}{c}
      \rho       \\
      \rho u
    \end{array} \right] \;, \\[10pt]
   { \bf F} \equiv  
    \left[ \begin{array}{c}
      f_1       \\
      f_2
    \end{array} \right]
     \equiv 
    \left[ \begin{array}{c}
      \rho u               \\
      \rho u^2 + a^2 \rho
    \end{array} \right] \;. 
  \end{array} \right\}   
\end{equation}
Here $\rho$ is density, $u$ is velocity and $a$ is a prescribed constant speed of sound.   The Jacobian, eigenvalues and corresponding right  eigenvectors are
\begin{equation}                                      \label{eq:Roe14}  
      \left.\begin{array}{c} 
      \bf A(\bf Q) =
      \left[ \begin{array}{cc}
       0          & 1     \\
      a^2 - u^2  & 2u
      \end{array} \right] \;, \\[8pt]
      \lambda_{1}=u-a \;, \hspace{3mm} \lambda _{2}=u+a \;, \\[8pt]                
      {\bf K}^{(1)}=
      \left[
      \begin{array}{c}
      1 \\
      u-a
      \end{array}\right] \;, \hspace{3mm} 
      {\bf K}^{(2)}=
      \left[
      \begin{array}{c}
      1 \\
      u+a
     \end{array}\right]   \;.
     \end{array}\right\}
\end{equation}
Following Roe \cite{Roe:1981a}  the {\em parameter vector} is chosen as
\begin{equation}                                              \label{eq:Roe15}
  \bf W \equiv
  \left[ \begin{array}{c}
    w_1 \\ w_2
  \end{array} \right]
  \equiv \frac {\bf Q}{\sqrt \rho} =\sqrt{\rho}
  \left[ \begin{array}{c}
    1 \\  u
  \end{array} \right] \;.
\end{equation}
It is easily verified that  $\bf Q$ and $\bf F$ can be expressed in terms of the components
$w_1, w_2$ of $\bf W$, namely
\begin{equation}                                              \label{eq:Roe16}
  {\bf Q }\equiv
  \left[ \begin{array}{c}
    q_1 \\ q_2
  \end{array} \right]
  \equiv  w_1 {\bf W} =
  \left[ \begin{array}{c}
    w_1^2 \\ w_1 w_2
  \end{array} \right]
\end{equation}
and
\begin{equation}                                              \label{eq:Roe17}
  {\bf  F} \equiv
  \left[ \begin{array}{c}
    f_1 \\ f_2
  \end{array} \right]
  \equiv  
  \left[ \begin{array}{c}
    w_1 w_2 \\ w_2^2 + a^2w_1^2
  \end{array} \right] \;.
\end{equation}
One now looks for an averaged parameter vector $\tilde{\bf W} = [\tilde w_1, \tilde w_2]^T$.  It may seem surprising but this is found by  simple {\em arithmetic averaging}
\begin{equation}                                              \label{eq:Roe18}
  \tilde{\bf W} =
  \left[ \begin{array}{c}
    \tilde w_1 \\ \tilde w_2
  \end{array} \right]
  = \frac 1 2 ({\bf W}_{L} + {\bf W}_{R} )\;.
  \end{equation}
In fact
\begin{equation}                                              \label{eq:Roe19}
   \tilde{\bf W} =
  \left[ \begin{array}{c}
    \tilde w_1 \\ \tilde w_2
  \end{array} \right]= \frac 1 2
  \left[ \begin{array}{c}
    \sqrt{ \rho_L } + \sqrt{ \rho_R } \\ 
    \sqrt{ \rho_L } u_L + \sqrt{\rho_R} u_R  
  \end{array} \right] \;.
\end{equation}
Then,  two matrices $\tilde{\bf B} = \tilde{\bf B}(\bf{\tilde W})$ and $\tilde{\bf C}
 = \tilde{\bf C}(\bf {\tilde W})$  are found, such that the jumps $\Delta \bf Q$ and  $\Delta \bf  F$ in (\ref{eq:Roe12}) can be expressed in terms of the  jump $\Delta \bf W$, namely
\begin{equation}                                              \label{eq:Roe20}
  \Delta {\bf Q} = \tilde{\bf B} \Delta {\bf W} \;, ~~
  \Delta {\bf  F} = \tilde{\bf C} \Delta {\bf W} \;.
\end{equation}
Use of these two expressions produces
\begin{equation}                                              \label{eq:Roe21}
  \Delta {\bf F} = (\tilde{\bf C} \tilde{\bf B} ^ {-1} ) \Delta{\bf Q} \;,
\end{equation}
which if compared with condition (\ref{eq:Roe10}) produces the sought Roe averaged matrix
\begin{equation}                                              \label{eq:Roe22}
  \tilde{\bf A} = \tilde{\bf C} \tilde{\bf B}^{-1} \;.
\end{equation}
Matrices $\tilde{\bf B}$ and $\tilde{\bf C}$ satisfying (\ref{eq:Roe20}) are 
\begin{equation}                                              \label{eq:Roe23}
  \tilde{\bf B} =
  \left[ \begin{array}{cc}
    2 \tilde w_1   & 0      \\
    \tilde w_2     & \tilde w_1
  \end{array} \right]  \;, ~
  \tilde{\bf C} =
  \left[ \begin{array}{cc}
    \tilde w_2          & \tilde w_1      \\
    2a^2 \tilde w_1     & 2 \tilde w_2^2
  \end{array} \right]  \;,
\end{equation}
which the reader can easily verify. The sought Roe matrix is then
\begin{equation}                                              \label{eq:Roe24}
  \tilde{\bf A} =
  \left[ \begin{array}{cc}
    0                    & 1                       \\
    a^2 - \tilde u^2     & 2 \tilde u
  \end{array} \right]  \;,
\end{equation}
where $\tilde u$ is the {\bf Roe averaged velocity}, given as
\begin{equation}                                              \label{eq:Roe25}
  \tilde u = 
  \frac{
    \sqrt{\rho_L} u_L + \sqrt{\rho_R} u_R
  }{
    \sqrt{\rho_L} + \sqrt{\rho_R}
  }     \;.
\end{equation}
Compare (\ref{eq:Roe24}) with the matrix in (\ref{eq:Roe14}).  As the sound speed $a$ is constant, no averaged  $\tilde \rho$ is required. 

Having found $\tilde{\bf A}$ one computes the averaged eigenvalues and eigenvectors 
\begin{equation}                                              \label{eq:Roe26}
  \tilde \lambda_1 = \tilde u - a \;, ~~
  \tilde \lambda_2 = \tilde u + a
\end{equation}
and
\begin{equation}                                              \label{eq:Roe27}
  \tilde{\bf K}^{(1)} =
  \left[ \begin{array}{c}
    1 \\ \tilde u - a  
  \end{array} \right] \;, ~~
  \tilde{\bf K}^{(2)} =
  \left[ \begin{array}{c}
    1 \\ \tilde u + a  
  \end{array} \right] \;.
\end{equation}
The eigenvalues are real and distinct and hence the right eigenvectors are linearly independent. 

As seen in \cite{Toro:2009a} the numerical flux for a linear system of $m$ equations
with eigenvalues $\tilde \lambda_i $, wave strengths  $\tilde \alpha_i $ and eigenvectors  $\tilde{\bf K}^{(i)}$ can be expressed as
\begin{equation}                                              \label{eq:Roe28}
  {\bf F}_{i+ \frac 1 2} = {\bf F}_L 
  + \sum_{\tilde{\lambda}_i \leq 0} \tilde \alpha_i \tilde \lambda_i \tilde{\bf K}^{(i)} \;, 
\end{equation}
or
\begin{equation}                                              \label{eq:Roe29}
 {\bf F}_{i+ \frac 1 2} = {\bf F}_R 
  - \sum_{\tilde{\lambda}_i \geq 0} \tilde \alpha_i \tilde \lambda_i \tilde{\bf K}^{(i)} \;.
\end{equation}
Alternatively,  it may also be written thus
\begin{equation}                                              \label{eq:Roe30}
  {\bf F}_{i+ \frac 1 2} = \frac 1 2 ({\bf F}_L + {\bf F}_R ) 
  - \frac 1 2 
  \sum_{i = 1}^m \tilde \alpha_i | \tilde \lambda_i | \tilde{\bf K}^{(i)} \;.
\end{equation}
For our example,  the wave strengths  $\tilde{\alpha}_i$ are given by the $2 \times 2$ linear system
\begin{equation}                                              \label{eq:Roe31}
  \Delta {\bf Q} \equiv
  \left[ \begin{array}{c}
    \Delta q_1 \\ \Delta q_2
  \end{array} \right]
  = \sum_{i= 1}^2 \tilde \alpha_i \tilde{\bf K}^{(i)} \;.
\end{equation}
The solution is easily verified to be
\begin{equation}                                              \label{eq:Roe32}
  \left. \begin{array}{ccc}
    \tilde \alpha_1 & = & 
    \displaystyle{
      \frac{
        \Delta q_1 ( \tilde u + a ) - \Delta q_2
      }{
        2a
      }   
    } \;,                 \\[3mm]
    \tilde \alpha_2 & = & 
    \displaystyle{
      \frac{
        - \Delta q_1 ( \tilde u - a ) + \Delta q_2
      }{
        2a
      }   
    }\;,                 
  \end{array} \right\}
\end{equation}
with the obvious definitions
\begin{equation}                                              \label{eq:Roe33}
 \Delta q_1 \equiv \rho_R - \rho_L, ~~ \Delta q_2 \equiv \rho_R u_R - \rho_L u_L \;.
\end{equation}
The corresponding Roe numerical flux $\bf F_{i+\frac 1 2}$ now follows from 
using (\ref{eq:Roe25})--(\ref{eq:Roe27}) and (\ref{eq:Roe32})--(\ref{eq:Roe33})
into any of the expressions (\ref{eq:Roe28})--(\ref{eq:Roe30}).

Details for the three-dimensional Euler equations can be found in the original Roe paper \cite{Roe:1981a} and in Chap. 11 of \cite{Toro:2009a}.

\subsection{ Closing remarks on Roe's solver}

\begin{enumerate}

\item There is an alternative formulation of the method due to Roe and Pike.  For full details see \cite{Roe:1984a}.

\item The numerical scheme (\ref{eq:Roe2}) with numerical flux (\ref{eq:Roe28}), (\ref{eq:Roe29}) or (\ref{eq:Roe30}) is a first-order accurate upwind scheme.   For many standard test problems, Roe's approximate solver produces results closely resembling those obtained using an exact Riemann solver, although significant differences can arise in certain flow regimes.

\item In \cite{Roe:1981a} Roe already noted a limitation of his linearised Riemann solver, namely the entropy property.  In the presence of transonic flow, or transcritical flow,  the solver as stated, instead of computing a transonic rarefaction wave,  it computes an entropy violating (unphysical) shock.   
Harten and Hyman \cite{Harten:1983c},   Roe and Pike  \cite{Roe:1984a}, Roe \cite{Roe:1992b}, Dubois and Mehlman \cite{Dubois:1993a} and others, have produced  appropriate modifications to the scheme to overcome this limitation.   Details on how to construct an {\it entropy fix} for the Euler equations  are found in Chapter 11 of \cite{Toro:2009a}. 

\item Another issue concerns the difficulty of {\bf incipient vacuum}.  In \cite{Einfeldt:1991a},  corrections to the basic  Roe scheme to avoid this problem were presented; the authors also  showed that in fact this anomaly afflicts {\em all linearised Riemann solvers}. 
\end{enumerate}

\section{Nonlinear schemes: limiters}
\label{sec:NonlinearSchemes}

Already in his 1981  paper \cite{Roe:1981a} Roe expressed the desirability of increasing the accuracy of the first-order schemes he had proposed to solve hyperbolic equations.   He was aware of the major theoretical barrier stated in Godunov's theorem \cite{Godunov:1959a}    
concerning the impossibility of constructing {\bf linear} monotone schemes  of accuracy greater than one.  So the task was to construct {\bf nonlinear} schemes of second order, or better still, third order.  

In this section we briefly review some of Roe's main contributions to the construction of nonlinear schemes.

\subsection{Godunov's theorem}

Consider the model PDE for the unknown $q(x,t)$ written in conservative form
\begin{equation}                                              \label{eq:NLRoe1}
 \partial_{t}q+ \partial_{x}f(q)=0 \;,  \hspace{2mm} f(q)=\lambda q \;, 
\end{equation}
where $\lambda$ is a prescribed constant wavespeed.  Consider also the class of 
{\bf linear schemes}
\begin{equation}                                                  \label{eq:NLRoe2}
    q_{i}^{n+1} = \sum_{k=-l}^{k=r}b_{k}q_{i+k}^{n} \;
\end{equation}
%{}
to solve (\ref{eq:NLRoe1}),  with constant coefficients $b_{k}$, where $l$ and $r$ are non-negative integers.    A {\bf linear scheme}  (\ref{eq:NLRoe2}) is said to be {\bf monotone} if the coefficients are non-negative,  that is $b_{k}\ge 0$.   

Recall from (\ref{eq:Roe2}) the class of conservative schemes, which for  (\ref {eq:NLRoe1}) take the form
\begin{equation}                                                  \label{eq:NLRoe3}
 q_{i}^{n+1}=q_{i}^{n} -\frac{\Delta t}{\Delta x}(f_{i+\frac{1}{2}}-f_{i-\frac{1}{2}})\;, 
\end{equation}
where $f_{i+\frac{1}{2}}$ is the numerical flux defining the scheme.   In the Godunov method the numerical flux for (\ref{eq:NLRoe1}) is
\begin{equation}                                                  \label{eq:NLRoe3a}
 f_{i+\frac{1}{2}}=\left\{\begin{array}{ccc}
 \lambda q_{i}^{n}       & \mbox{if} & \lambda \ge 0 \;, \\
 \lambda q_{i+1}^{n}  & \mbox{if}  & \lambda \le 0 \;. 
 \end{array}\right.
\end{equation}
This results from solving the appropriate local Riemann problem.  For background see \cite{Toro:2009a}.
A special case of (\ref{eq:NLRoe2}) is the family of three--point linear schemes
\begin{equation}                                                  \label{eq:NLRoe3b}
    q_{i}^{n+1} = b_{-1} q_{i-1}^{n} + b_{0} q_{i}^{n} + b_{1} q_{i+1}^{n} \;. 
\end{equation}
For the first-order Godunov upwind scheme written as (\ref{eq:NLRoe3b}),   there are two cases.
\begin{enumerate}
\item For $\lambda>0$: $b_{-1}=c$,  $b_{0}=1-c$ and $b_{1}=0$.  
\item For $\lambda<0$: $b_{-1}=0$, $b_{0}=1+c$ and $b_{1}=-c$.   
\end{enumerate}
Here 
\begin{equation}                                                  \label{eq:NLRoe3c}
    c=\frac{\lambda \Delta t}{\Delta x}
\end{equation}
is the Courant number.  Another example of (\ref{eq:NLRoe3b}) is the  Lax-Wendroff method, in which case
$b_{-1}=\frac{1}{2}c(1+c)$, $b_{0}=(1-c^{2})$ and  $b_{1}=-\frac{1}{2}c(1-c)$.  

The reader is encouraged to verify that the Godunov method is always monotone under the condition $|c|\le 1$, which is also its linear stability condition. On the other hand, the Lax-Wendroff method is never monotone and will therefore produce spurious oscillations outside the range $[q_{min}, q_{max}]$,  containing the initial conditions.   Exceptional cases are $c=0$ and $|c|=1$. \\

\noindent{\bf Remark:} There are numerical schemes  defined as monotone schemes in the literature nevertheless produce {\bf internal spurious oscillations},   that is, the solution may oscillate while still remaining within the interval $[q_{min}, q_{max}]$. The Lax-Friedrichs method is a known example, though not very well-known, in the literature 
 \cite{Toro:2025a}.\\

\noindent{\bf  An accuracy lemma.} A linear scheme (\ref{eq:NLRoe2}) is said to be $p$-th order accurate in space and time ($p \geq 0$), in the sense of local truncation error, if and only if
\begin{equation}                                                  \label{eq:NLRoe4}
      \sum_{k=-l}^{r} k^{\eta} b_k = (-c)^{\eta} \;, \hspace{3mm} \eta = 0, 1,\ldots, p \;. 
\end{equation} 
Note that for purely notational convenience one introduces the notation $0^{0}=1$ (mathematically meaningless).\\

\noindent {\bf Proof.}  For details see Roe \cite{Roe:1981c} and Toro \cite{Toro:2009a},  for example.  
For an extension of the accuracy lemma to two and three space dimensions see Billett and Toro \cite{Billett:1997a} or  Toro \cite{Toro:2009a}.\\

\noindent{\bf Godunov's Theorem.} There are no monotone, linear schemes (\ref{eq:NLRoe2}) for the linear advection equation in (\ref{eq:NLRoe1}) with constant $\lambda$, of accuracy two or higher
 \cite{Godunov:1959a}.\\

\noindent {\bf Proof.}
For details on the proof see Roe \cite{Roe:1981c}; see also Toro \cite{Toro:2009a}.

\subsection{Nonlinear schemes \'a la Roe}

One of the first reported attempts at providing a rational approach to circumventing Godunov's theorem is due to  Roe \cite{Roe:1981c}. The central idea is to construct adaptive algorithms that would adjust themselves to the local 
nature of the solution. This results, obviously, in variable coefficient (non--linear) schemes, even when applied to  linear PDEs. The coefficients must be functions of the data.  The following definition relates an algorithm for the linear equation (\ref{eq:NLRoe1}) to particular classes of data.\\

\noindent{\bf Data compatible scheme.}
A scheme is compatible with a given data set $\{q_i^n\}$ if the solution 
$q_i^{n+1}$ at each point $i$, as given by the algorithm, is bounded by the
{\bf upwind pair} $(q_{i-s}^n, q_i^n)$, where $s \equiv \mbox{sign}(c) = \mbox{sign}(\lambda)$.

Note that for a given point $i$ the upwind pair is $(q_{i-1}^n, q_i^n)$ for $c>0$ and $(q_i^n, q_{i+1}^n)$ for $c<0$. The characteristic through $(i, n+1)$ of the PDE
(\ref{eq:NLRoe1}) intersects the real line at a point that lies between the two data points of the upwind pair. One may therefore rephrase the above definition as follows: a scheme is compatible with the data if $q_i^{n+1}$ lies between $q_{i-s}^n$ and $q_i^n$. 

The data compatibility requirement  of an algorithm is conveniently expressed as follows:
\begin{equation}                                            \label{eq:NLRoe8}
  \min\{ q_{i-s}^n, q_i^n \} \leq q_i^{n+1} \leq \max\{q_{i-s}^n, q_i^n\} \;,
\end{equation}
which in fact is equivalent to requiring
\begin{equation}                                            \label{eq:NLRoe9}
  0 \leq \frac{ q_i^{n+1}- q_i^n }{ q_{i-s}^n-q_i^n } \leq 1 \;.
\end{equation}
This can be seen as follows:  assume $q_{i-s}^n \leq q_i^n$.  Imposing (\ref{eq:NLRoe8}) gives
$q_{i-s}^n \leq q_i^{n+1} \leq q_i^n$, that is $q_{i-s}^n - q_i^n \leq q_i^{n+1} - q_i^n \leq 0$.
Division through by $q_{i-s}^n-q_i^n < 0$ gives (\ref{eq:NLRoe9}). The case
for which $q_{i-s}^n \geq q_i^n$ follows similarly.\\

We note that Godunov's flux  for (\ref{eq:NLRoe1}) may also be written as
\begin{equation}                                            \label{eq:NLRoe10}
  f_{i+\frac 1 2} = \frac 1 2 (1+s)(\lambda q_i^n) + \frac 1 2 (1-s)(\lambda q_{i+1}^n)  \;.
\end{equation}
Then, the conservative scheme (\ref{eq:NLRoe3})  gives
\begin{equation}                                            \label{eq:NLRoe11}
  q_i^{n+1} = q_i^n + sc(q_{i-s}^n - q_i^n ) \;.
\end{equation}
Note that the data compatibility condition (\ref{eq:NLRoe9}) for (\ref{eq:NLRoe11}) may conveniently be written as
\begin{equation}                                            \label{eq:NLRoe12}
  0 \leq \frac{ q_i^{n+1} - q_i^n }{ q_{i-s}^n - q_i^n } \le |c| \leq 1 \;.
\end{equation}
That is, the Godunov first--order upwind method is compatible with all possible sets of data. This is not a property enjoyed by the Lax--Wendroff method.  See \cite{Toro:2009a} for details.

Roe \cite{Roe:1981c} observed that given the set $S_{DC}$ of all possible data states one can select combinations of schemes such that the subsets for which the schemes satisfy the data compatibility requirement do cover the full
set $S_{DC}$; in this way he constructed {\em adaptive, non--linear schemes that are monotone and second  order accurate},  thus circumventing Godunov's theorem.
More recent applications of the data compatibility concept impose directly the data compatibility requirement (\ref{eq:NLRoe12}) on schemes that depend on some free parameter to allow for variable coefficients. 

\subsection{Flux limiters}

One way to understand Roe's contributions to the construction of nonlinear methods is to adopt the {\em flux limiter approach}; see Roe  \cite{Roe:1983a} and Sweby \cite{Sweby:1984a}. 
It is worth noting that, during his time at RAE, Phil Roe served as the industrial PhD supervisor of Peter Sweby at the University of Reading.

The flux limiter approach requires a high-order flux $f_{i+\frac 1 2}^{HI}$ associated with a scheme of accuracy greater than or equal to two and a  low-order flux $f_{i+\frac 1 2}^{LO}$ associated with a {\em monotone},  first--order scheme.   The reader is warned that not all first-order accurate methods are monotone.  It is also noted that the Godunov first-order upwind method  is the most accurate scheme in the family of all monotone methods; it has the smallest local truncation error.

We assume the model equation (\ref{eq:NLRoe1}) as solved by the conservative scheme (\ref{eq:NLRoe3}).  One then defines a high--order  {\bf limited} flux as
\begin{equation}                                          \label{eq:NLRoe20}
  f_{i+ \frac 1 2}^{Lim}
  = f_{i+ \frac 1 2}^{LO} 
  + \psi_{i + \frac 1 2}
  [ f_{i+\frac 1 2}^{HI} - f_{i+\frac 1 2}^{LO}] \;,
\end{equation}
where $\psi_{i+\frac 1 2}$ is a {\bf flux limiter} function yet to be determined. 
Three crucial steps are involved: 
\begin{enumerate}
\item A selection of specific fluxes $f_{i+\frac 1 2}^{LO}$ and $f_{i+\frac 1 2}^{HI}$ must be made,
\item Substitution of limited fluxes  $f_{i-\frac 1 2}^{Lim}$ and $f_{i+\frac 1 2}^{Lim}$ into the conservative formula (\ref{eq:NLRoe3}),
\item  Imposition of the data-compatibility condition (\ref{eq:NLRoe9}).
\end{enumerate}
Then, the algebraic task at hand is to determine a {\bf TVD} (Total Variation Diminishing) region for the limiter $\psi(r)$, where the parameter $r$ monitors the variation of the data and is given as
\begin{equation}                                            \label{eq:NLRoe25}
  r = \frac{ \Delta_{upw} }{ \Delta_{loc} } = 
  \left\{ \begin{array}{ll}
    \displaystyle{
      \frac{ q_i^n - q_{i-1}^n }{ q_{i+1}^n - q_i^n }     
    }  
    \;, & \lambda > 0                        \;,    \\[3mm]
    \displaystyle{
      \frac{ q_{i+2}^n - q_{i+1}^n }{ q_{i+1}^n - q_i^n } 
    }
    \;, & \lambda < 0      \;.
    \end{array} \right.
\end{equation}
In general,  $\psi_{B} \le \psi(r) \le \psi_{T}$ and the complete determination of $\psi(r)$ depends on specific choices for the bottom $\psi_{B}$ and top $\psi_{T}$ bounds.  
\begin{figure}
      \centerline{\includegraphics[width=7cm,angle=0]{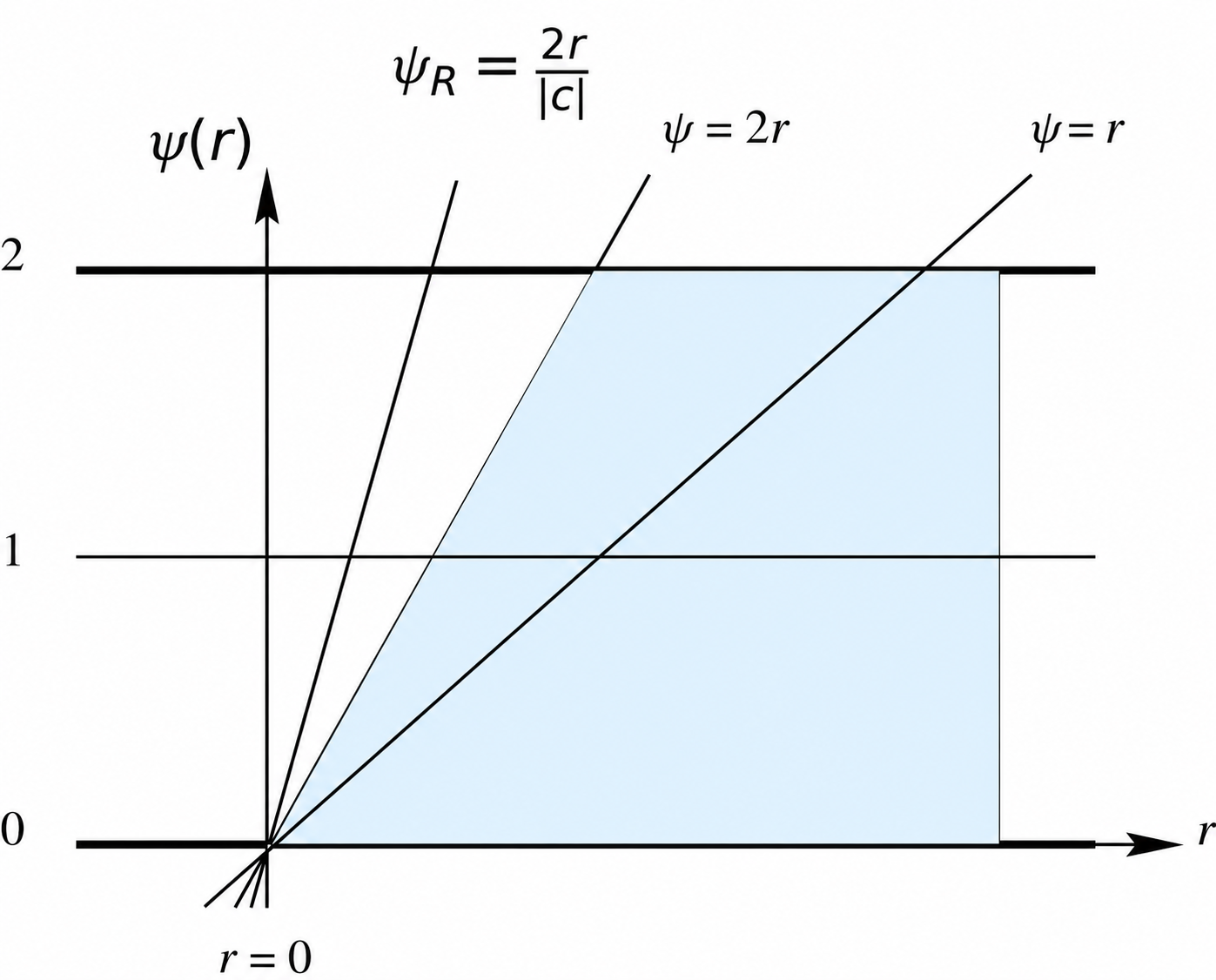}}
      \caption{Sweby's  TVD Region (light blue) \cite{Sweby:1984a}.  The line $\psi_{R}$, which depends on 
      $r$ and $c$,  gives rise to a larger TVD region not discussed here.   
      The bottom and top boundaries are $\psi=0$ and $\psi=2$ respectively,  so that the 
      Sweby TVD region and derived  limiters are  independent of Courant number $c$.}
      \label{fig:SwebyTVD}
\end{figure}

The classical Sweby TVD region  \cite{Sweby:1984a} is depicted in  Fig.   \ref{fig:SwebyTVD} (in blue).
This is a special case of larger TVD regions derived in Chap.  13 of \cite{Toro:2009a} and results from choosing the bounds
\begin{equation}                                            \label{eq:NLRoe32}
    \psi_T = 2 \leq \frac 2 {1-|c|}  \;, \hspace{3mm} \forall c \;; \hspace{3mm} \psi_B = 0 \;.
\end{equation}
There are three horizontal lines: the bottom and top boundaries and the second-order line $\psi=1$.  There are three inclined lines of positive slopes $2r/|c|$, $2r$ and $r$. Note that $\psi_{R}=2r/|c|$ lies outside the Sweby TVD region. 

\begin{figure}
      \centerline{\includegraphics[width=7.5cm,angle=0]{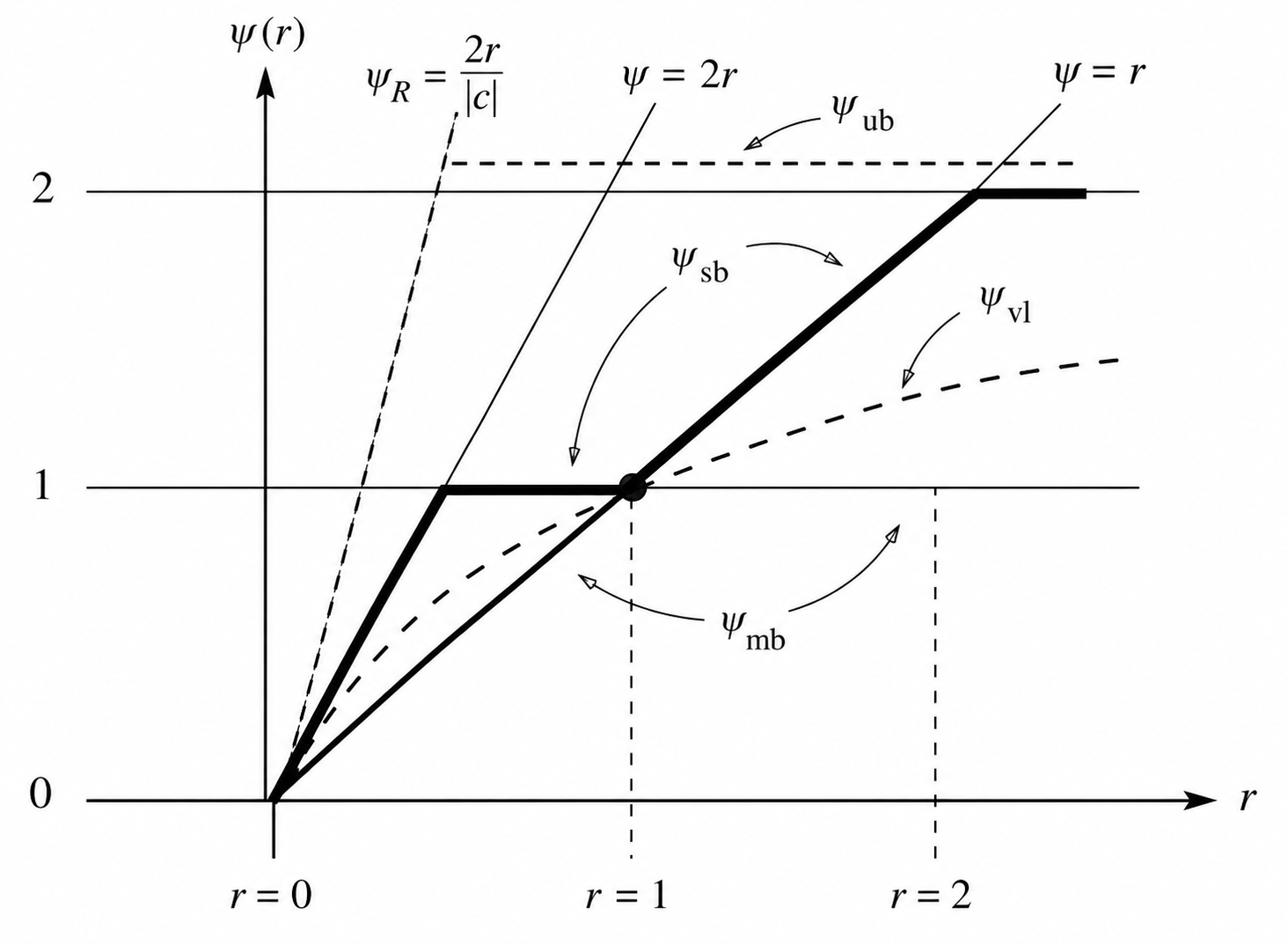}}
      \caption{Four flux limiter functions constructed from the Sweby TVD region of Fig.   
      \ref{fig:SwebyTVD}. $\psi_{ub}$: ULTRABEE; $\psi_{sb}$: 
      SUPERBEE; $\psi_{vl}$: VanLeer and $\psi_{mb}$: MINBEE.}
      \label{fig:limiters}
\end{figure}

Two well-known limiter functions due to Roe and derived from the Sweby TVD region of Fig.  \ref{fig:SwebyTVD} are SUPERBEE
\begin{equation}                                            \label{eq:NLRoe33}
  \psi_{sb}(r) = 
  \left\{ \begin{array}{ll}
    0  \;,  & r \leq 0                      \;,   \\[2pt]
    2r \;,  & 0 \leq r \leq \frac 1 2       \;,   \\[2pt]
    1  \;,  & \frac 1 2 \leq r \leq 1       \;,   \\[2pt]
    r  \;,  & 1 \leq r \leq 2               \;,   \\[2pt]
    2  \;,  & r \geq 2                      \;, 
  \end{array} \right.
\end{equation}
and MINBEE
\begin{equation}                                            \label{eq:NLRoe34}
  \psi_{mb}(r) = 
  \left\{ \begin{array}{ll}
    0 \;, & r \leq 0                 \;,        \\[2pt]
    r \;, & 0 \leq r \leq 1          \;,      \\[2pt]
    1 \;, & r \geq 1                 \;. 
  \end{array} \right.
\end{equation}
Fig.  \ref{fig:limiters}  illustrates four flux limiters constructed from Sweby's  TVD region of Fig.  \ref{fig:SwebyTVD}, namely $\psi_{ub}$: ULTRABEE; $\psi_{sb}$: SUPERBEE; $\psi_{vl}$: VanLeer and $\psi_{mb}$: MINBEE.

\subsection{Sample numerical results}

\begin{figure}
\centerline{
\includegraphics[height=75mm,angle=0]{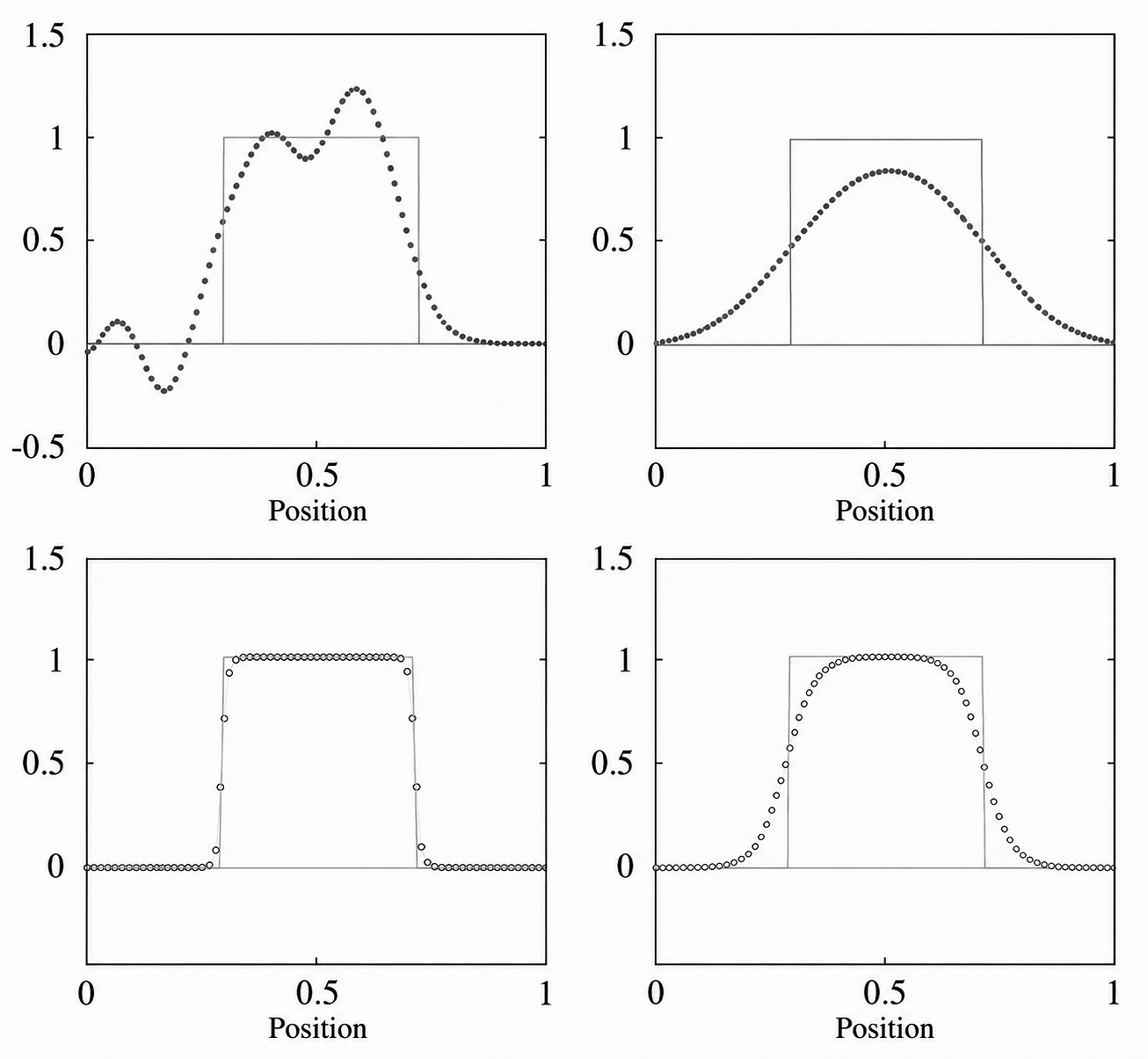}
}
\vspace{0mm}
\caption{Numerical (symbols) and exact (line) solutions for square wave test 
         (\ref{eq:NLRoe35}) at time $t=10$. 
        Top left: Lax-Wendroff; top right:  Godunov first-order upwind; bottom left: SUPERBEE; 
        bottom right: MINBEE.    Mesh used M=100 cells; $C_{cfl}=0.8$
         }
         \label{fig:LAresults}
\end{figure}

We consider a test problem  for the linear advection (\ref{eq:NLRoe1}) with $\lambda=1$, initial condition (\ref{eq:NLRoe35}) and periodic boundary conditions.  Computations are carried out on mesh $M=100$, with CFL coefficient $C_{cfl}=0.8$ and output time $t=10.0$. 
\begin{equation}                                                  \label{eq:NLRoe35}
       q(x,0)=
       \left\{\begin{array}{rcl}
       0              &     \mbox{ if } &0 \le  x \leq 0.3 \;,\\
       1              &     \mbox{ if } &  0.3 < x < 0.7 \;,\\
       0              &     \mbox{ if } &  0.7 \le x \le 1.0  \;. 
       \end{array} \right.
\end{equation}
In Fig.  \ref{fig:LAresults} we compare the numerical solution (symbols) with the exact solution (full line).  The top left frame shows results for the Lax-Wendroff method (linear, non-monotone). The top right frame shows results from the  Godunov first-order upwind method (monotone). The bottom left frame shows results from the TVD method with the SUPERBEE  limiter and the bottom right frame shows results from the TVD method with the  MINBEE limiter. Obviously, the results from the TVD methods,  bottom frames,  are superior to those above.  That is the achievement of flux limiter methods. 

\subsection{Remarks on nonlinear methods}

Research on nonlinear methods to circumvent Godunov's theorem \cite{Godunov:1959a}, during  the 1970s and 1980s,  is substantial.  Kolgan,  in 1972,   pioneered these developments \cite{Kolgan:1972a}.   Harten and others developed the theory of Total Variation Diminishing (TVD) methods \cite{Harten:1983a},\cite{Harten:1984a}.  Information on further contributions on related topics is found in Chap. 13 of \cite{Toro:2009a}, for example.  Back in 1989, Sweby \cite{Sweby:1989a} widened the discussion on TVD methods for equations containing source terms.

The extension to nonlinear systems is empirical but results for practical problems are quite satisfactory.  As these methods are at most second-order accurate,  their performance for wave propagation problems over long times and distances,  like in acoustics,  is not very satisfactory.  For shock-dominated problems these methods remain competitive. 

More recently, there have been remarkable developments in very-high-order numerical methods,  well above two; see the recent textbook \cite{Toro:2025a} and the many references therein.  The nonlinear versions of these methods rely on nonlinear spatial reconstructions.  But  interest in the {\bf old limiters} has resurfaced in recent years, also in the context of schemes of very high order of accuracy.

\section{Upwinding Source Terms}

Phil Roe seems to have been the first researcher to consider the problem of extending the upwind philosophy to the discretisation of hyperbolic equations with source terms \cite{Roe:1986c}. In this section, I first briefly review Roe's basic ideas on the subject using a scalar model equation and relate them to the subsequent early work of Berm\'udez and V\'azquez \cite{Bermudez:1994a}.  I then discuss some alternative ways of viewing the subject and, finally, present computational results obtained from applying Roe's ideas to a nonlinear system with source terms.

\subsection{Roe's early ideas}

In his seminal paper of 1986  \cite{Roe:1986c} Roe considered the following model equation for the unknown function $q(x,t)$
\begin{equation}                          \label{eq:RoeS1}
    \partial_{t}q(x,t)+\lambda \partial_{x}q(x,t)=s\;,  \hspace{3mm} s =g(x) \;,
\end{equation}
for a constant,  positive wavespeed $\lambda$. The source term $s$ is simply an {\it external} function $g(x)$ of distance $x$ alone; that is the source term is independent of the unknown $q(x,t)$ of the differential equation.  This particular choice for the source term could reasonably be related to compressible flow in a tube with area variation \cite{Toro:2009a} and to shallow water flows with variable bathymetry \cite{Toro:2024a}.
Such a relation becomes more evident for particular forms of the governing equations.
 
The notion of upwinding is immediately revealed from the exact solution of (\ref{eq:RoeS1}), namely
\begin{equation}                          \label{eq:RoeS2}
      q(x,t) =  q^{(0)}(x - \lambda t) + \frac{1}{\lambda} \int_{x - \lambda t}^{x} g(\xi) \; {\rm d} \xi \;,
\end{equation}
where $q^{(0)}(x)$ is the initial condition.  In a finite difference setting Roe put in evidence the  upwind flavour of solution (\ref{eq:RoeS2}) through the integration range $[x - \lambda t,x]$ biased in the upwind direction that affects both the advection term of the equation as well as the source term.

Now we depart slightly from Roe's original example and assume a source term of the form 
\begin{equation}                          \label{eq:RoeS3} 
    s(q)=\beta q \;,  
\end{equation}
with $\beta$ a constant, satisfying  $\beta \le 0$. One then may conveniently express the discrete version suggested by (\ref{eq:RoeS2}) as
\begin{equation}                          \label{eq:RoeS4} 
    q_{i}^{n+1}  =  q_{adv}  + \Delta t s(q_{i}^{s})  \;,
\end{equation}
where $q_{adv}$ represents a discretisation of the advection term in the PDE and the second term represents  the source term $s(q)$,  with $q_{i}^{s}$ a {\bf source state},  an approximation to $q$ intended to account for upwinding. 

Roe considered various options for evaluating the integral (\ref{eq:RoeS2}),  from which the following upwind discretisations for $q_{i}^{s}$ are suggested:
\begin{equation}                          \label{eq:RoeS5} 
     q_{i}^{s}= \omega q_{i-1}^{n} + (1-\omega) q_{i}^{n} \;, \hspace{2mm} 0\le \omega \le 1 \;.
\end{equation}
The upwind influence when $\lambda>0$ comes through the term $\omega q_{i-1}^{n}$; the role of the centred part $(1-\omega) q_{i}^{n}$ is reduced by the effect of $\omega$.
For  example,  by using the Godunov upwind scheme for the advection term $q_{adv}$ and the upwinded source term as suggested in (\ref{eq:RoeS5}), the complete scheme with upwinded advection and upwinded source term  reads
\begin{equation}                          \label{eq:RoeS6} 
\begin{array}{ccl}
    q_{i}^{n+1} & =&  c q_{i-1}^{n} +(1-c) q_{i}^{n}   \\
                           &  & + r \left[ \omega q_{i-1}^{n} + (1-\omega) q_{i}^{n} \right] \;,
\end{array}
\end{equation}
where $r=\Delta t \beta  \le 0$ is a non-dimensional {\it reaction number}.
Various choices of $\omega$ result in different upwind schemes.  Roe identified $\omega=\frac{1}{2}$ as a particularly interesting choice, as it leads to a second-order discretisation of the source term in the steady state.  Other choices considered by Roe will be discussed in the subsequent Secs.  \ref{sec:BV} and \ref{sec:Alter}.   In passing, we note that the particular choice $\omega=c$ in (\ref{eq:RoeS6}) reproduces identically the {\it splitting scheme} for source terms; see \cite{Toro:2024a},  page 291.

The first extensions of Roe's early ideas to nonlinear systems were reported by Glaister \cite{Glaister:1985a} and by Toro \cite{Toro:1986b}.  Illustrative computational results are given in Sec. \ref{sec:NLsystems}.

\subsection{Berm\'udez and V\'azquez early work}
\label{sec:BV}

In 1994,  Berm\'udez and V\'azquez \cite{Bermudez:1994a} published their study on schemes for discretizing source terms in hyperbolic balance laws.  They started  by considering the model equation
\begin{equation}                          \label{eq:BV1}
\partial_{t}q+ \partial_{x}f(q)=s(x,q(x,t))\;, \hspace{2mm} f(q)= \lambda q \;.
\end{equation}
They motivated the upwinding of the source term by integrating (\ref{eq:BV1}) along the  characteristic line passing through the point $(x_{i}, t_{n+1})$, for the case $\lambda>0$. They obtained
\begin{equation}                          \label{eq:BV2}
\begin{array}{c}
 q(x_{i},t_{n+1})  = q(x_{i}-\lambda\Delta t ,t_{n})+ \\[5pt]
 \displaystyle \int_{t_{n}}^{t_{n+1}}  \hspace{-7mm}s(x_{i}-\lambda(t_{n+1}-\tau),q(x_{i}-\lambda(t_{n+1}-\tau),\tau))d\tau \;.
  \end{array}
\end{equation}
The authors explored several formulae to evaluate the integral in (\ref{eq:BV2}) to obtain different expressions for upwinding the source term.  For the model equation with source term as in (\ref{eq:RoeS3}) their proposed schemes turn out to be identical to (\ref{eq:RoeS6}).  Beware of the sign difference in the source term (\ref{eq:RoeS3}) and that in the analysis of Berm\'udez and V\'azquez \cite{Bermudez:1994a}.
%
%\vspace{0.2cm}
\begin{table}[h] 
\begin{center}
\begin{tabular}{|c|c|c|} \hline
$\omega$ (Upwinding)      &   Lin.  Stab.    Cond.       \\ \hline \hline
$0$                                           &      $ c - \frac{1}{2} r \leq 1 $    \\ \hline
$\frac{1}{2} c $                      &       $ c  \leq 1 \;, \hspace{2mm}   r \geq -2 $  \\ \hline
$c$                                           &       $ c  \leq 1 \;, \hspace{2mm}   r \geq -2 $  \\ \hline
$\frac{1}{2}$                          &       $ c  \leq 1 \;, \hspace{2mm}   r \geq -2 $  \\ \hline
$1$                                           &       $ c  \leq 1 \;, \hspace{2mm}   r \geq -2 $  \\ \hline   
\end{tabular} 
\caption{Upwinding parameter $\omega$ for source term and corresponding linear stability conditions of resulting schemes.  Here $c$ is the Courant number and $r \le 0$ is the reaction number.}
%vspace{5mm}
\label{tab:BV}
\end{center}
\end{table}

Berm\'udez and V\'azquez \cite{Bermudez:1994a} also carried out a stability and positivity analysis of schemes (\ref{eq:RoeS6}). Their stability results are reproduced in Tab. \ref{tab:BV},  which reveals  that the centred scheme--no upwinding--has the most restrictive stability condition $c-\frac{1}{2}r \le 1$.  Upwinding the source term, that is $\omega>0$,  results in a more  generous stability condition, namely $c \le 1$ and $r \ge -2$.   The stability properties demonstrate that there is a net benefit in introducing upwinding in the discretisation of the source term.  See Chap. 13 of  \cite{Toro:2024a} for further discussion.

Berm\'udez and V\'azquez \cite{Bermudez:1994a} went further and proposed schemes for upwinding the source terms in nonlinear systems and demonstrated the benefits when  computing solutions to the shallow water equations with variable bathymetry.  See also the related papers \cite{Vazquez:1999a}, \cite{Toro:2002c}.

\subsection{Another way}
\label{sec:Alter}

Motivated by the early ideas presented above,  here I indulge in discussing some alternative  ways of upwinding the source term for the PDE
\begin{equation}                          \label{eq:T2}
\partial_{t}q+ \partial_{x}f(q)=s(q) \;.
\end{equation}
I adopt the finite volume framework and show that there are other possibilities for reproducing the schemes advanced in the early works of Roe \cite{Roe:1986c} and  Berm\'udez and V\'azquez \cite{Bermudez:1994a}.  Moreover, instead of recovering upwind information from tracing characteristics back in time, we do so by accounting for waves travelling forward in time via Riemann problem solutions.

We recall the finite volume approach.  As done for system (\ref{eq:Roe1}),  integration of (\ref{eq:T2})  in the control volume $V_{i}^{n}=[x_{i-\frac{1}{2}}, x_{i+\frac{1}{2}}] \times [t_{n}, t_{n+1}]$ gives
\begin{equation}                          \label{eq:T3}
    q_{i}^{n+1}=q_{i}^{n} - \frac{\Delta t}{\Delta x}(f_{i+\frac{1}{2}}-f_{i-\frac{1}{2}} ) 
    + \Delta t s_{i} \;,
\end{equation}
with definitions
\begin{equation}                           \label{eq:T4}
      \displaystyle{f_{i+\frac{1}{2}} = 
      \frac{1}{\Delta t} \int_{t_{n}}^{t_{n+1}} f(q(x_{i+\frac{1}{2}}, t) ) d t} \;
\end{equation}
and 
\begin{equation}                          \label{eq:T5}
      \displaystyle{s_{i} = \frac{1}{\Delta t \Delta x} \int_{t_{n}}^{t_{n+1}} \int_{x_{i-\frac{1}{2}}}^{x_{i+\frac{1}{2}}} s(q(x,t) ) d x d t} \;.
\end{equation}
As anticipated in Sec.  \ref{sec:NonlinearSchemes},  the classical upwind discretisation of the advection  term in (\ref{eq:T3}) results from numerical fluxes computed from the solution of the homogeneous Riemann problem
\begin{equation}                           \label{eq:T6}
\left.\begin{array}{l}
\partial_{t}q(x,t)+ \partial_{x}f(q(x,t))=0\\
 q(x,t_{n})=\left\{\begin{array}{lcc}
                q_{i}^{n}      & \mbox{if} & x < x_{i+\frac{1}{2}} \;, \\
               q_{i+1}^{n}  & \mbox{if}  & x > x_{i+\frac{1}{2}} \;. 
 \end{array}\right.
 \end{array}\right\}
\end{equation}
For the linear case $f(q)=\lambda q$,  the exact solution of  (\ref{eq:T6})  gives the Godunov upwind numerical flux as
\begin{equation}                            \label{eq:T7}
 f_{i+\frac{1}{2}}=\left\{\begin{array}{lcc}
 \lambda q_{i}^{n}       & \mbox{if} & \lambda \ge 0 \;, \\
 \lambda q_{i+1}^{n}  & \mbox{if}  & \lambda \le 0 \;. 
 \end{array}\right.
\end{equation}
Therefore, the scheme (\ref{eq:T3}) reads
\begin{equation}                          \label{eq:T8}
     q_{i}^{n+1}  =  q_{i}^{n} - c (q_{i}^{n} - q_{i-1}^{n}) +\Delta t s_{i} \;.
\end{equation}
The advection part is identical to that obtained from the finite difference interpretation of Godunov's method; see (\ref{eq:RoeS6}).  Defining the numerical source $s_{i}$ remains a pending task.   From here on we assume the model source term to be given by (\ref{eq:RoeS3}).

In order to evaluate $s_{i}$  from the integral (\ref{eq:T5}) we need to know the integrand $s(q(x,t))$ as a function of $x$ and $t$ inside the domain of integration.  And here there  is a problem.  In rigour we only know the solution at $t=t_{n}$, for which $q(x,t_{n})=q_{i}^{n}$ (cell average).  If we use this knowledge,  as a tentative approximation,  then the integral (\ref{eq:T5}) gives the {\bf numerical source} as
\begin{equation}                           \label{eq:T9}
s_{i}  = \beta q_{i}^{n}  \;.
\end{equation}
This type of source term discretisation has been called {\bf centred} discretisation. There is no upwinding,  no wave propagation information.

Let us look for a better approximation to the integrand in (\ref{eq:T5}).   We assume now that the integrand is approximated by the space-time dependent solution of the two neighbouring (homogeneous) Riemann problems  in the control volume $V_{i}^{n}$, respectively placed at  $x_{i-\frac{1}{2}}$ and $x_{i+\frac{1}{2}}$.  

\begin{figure}[h]
      \centerline{
            \includegraphics[scale=0.12, angle=0]{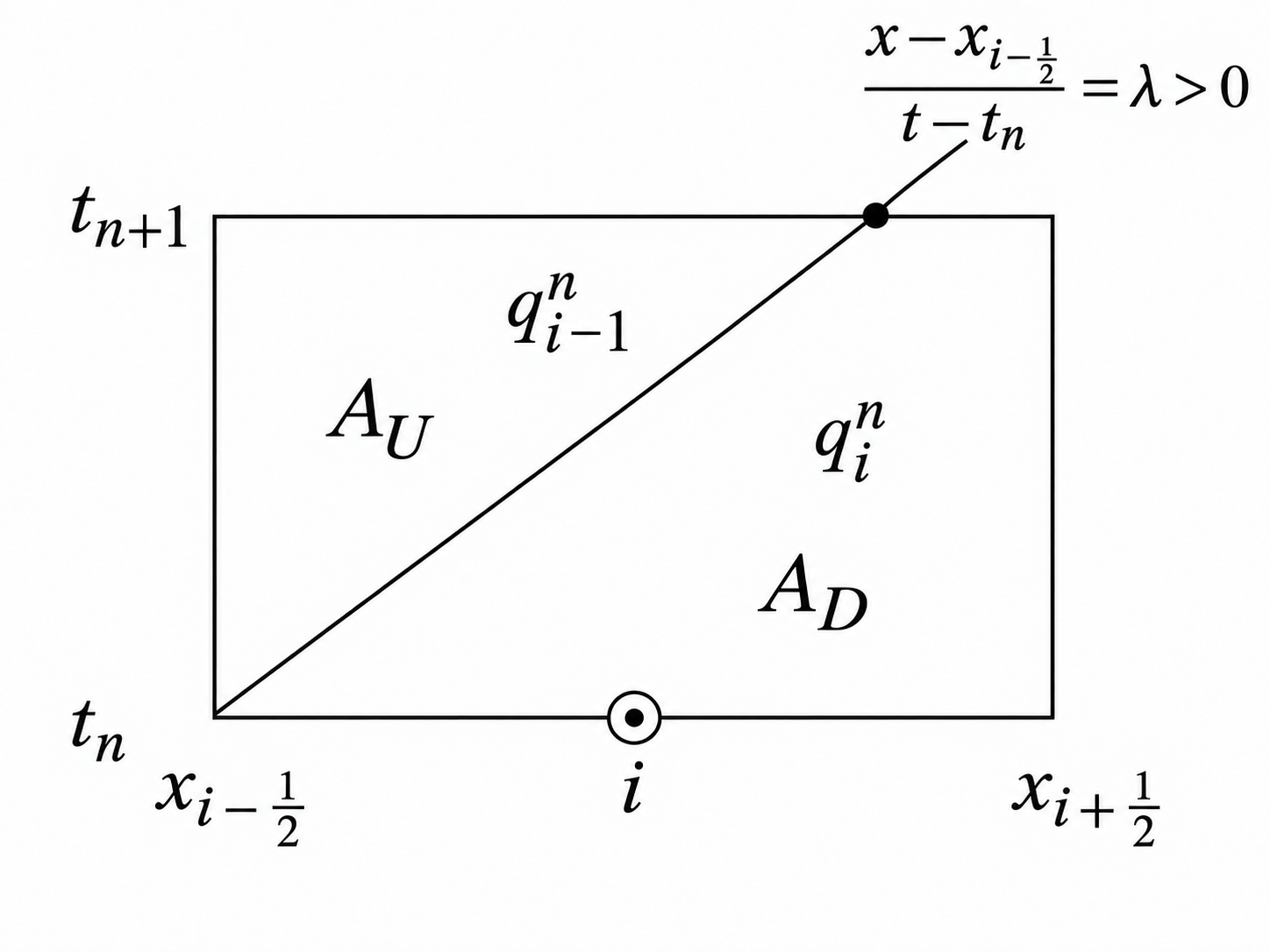}
      }
      \caption{Upwinding the source term.  For  $\lambda>0$ the solution of the left Riemann 
      problem at $x_{i-\frac{1}{2}}$ propagates the left data value $q_{i-1}^{n}$ to the right,  
      covering an area $A_{U}$ of the control volume.}
      \label{fig:source}
\end{figure}

Fig.   \ref{fig:source} illustrates the situation for  the linear case with $\lambda >0$,  for which only the left Riemann  problem  at $x_{i-\frac{1}{2}}$ affects the solution inside the control volume  $V_{i}^{n}$.  In fact the right Riemann problem at $x_{i+\frac{1}{2}}$ propagates $q_{i}^{n}$ into the right neighbouring volume with cell centre $x_{i+1}$,  and is therefore not relevant.   We remark that in practice, one works in local coordinates.

Within the control volume  $V_{i}^{n}$ there are two regions:  one of area $A_{U}$,  in which the solution takes on the value $q_{i-1}^{n}$,  and another with area $A_{D}$,  in which the solution takes on the value  $q_{i}^{n}$.   Straightforward evaluation of the integral (\ref{eq:T5}) gives
\begin{equation}                          \label{eq:T10}
 s_{i}  = \frac{\beta}{\Delta t \Delta x} [ A_{U} q_{i-1}^{n}+
 A_{D} q_{i}^{n} ] \;,
\end{equation}
which after simplifications reduces to
\begin{equation}                          \label{eq:T10b}
 s_{i}  =  \beta [ \frac{1}{2}cq_{i-1}^{n}+
(1-\frac{1}{2}c) q_{i}^{n}] \;.
\end{equation}
The first term in the bracket $\frac{1}{2}c q_{i-1}^{n}$ is the {\bf upwind contribution} due to wave propagation.  Note that (\ref{eq:T10b}) is identical to one of the schemes proposed by Roe and by Berm\'udez and V\'azquez.  See Table \ref{tab:BV}.

For  the case in which $\lambda <0$,  only the right Riemann problem at $x_{i+\frac{1}{2}}$ is relevant; its solution propagates the value $q_{i+1}^{n}$ inside the integration volume $V_{i}^{n}$; in the rest of $V_{i}^{n}$ the solution remains $q_{i}^{n}$.  In this case the source term becomes:
\begin{equation}                           \label{eq:T11}
 s_{i}  =  \beta[-\frac{1}{2}c q_{i}^{n}+(1+\frac{1}{2}c)q_{i+1}^{n}] \;.
\end{equation}

\paragraph{Fixed-time integrand.}
The approach just described can be generalised by aiming for an upwind source state $q_{i}^{s}$ so that the numerical source can be expressed as
\begin{equation}                            \label{eq:T12}
      s_{i}  =  s(q_{i}^{s})\;.
\end{equation}
Consider Fig.  \ref{fig:source} again, with the wave patterns from the two Riemann problems placed at $x_{i-\frac{1}{2}}$ and $x_{i+\frac{1}{2}}$.  Assume $q_{i}^{s}$ is computed as follows
\begin{equation}                           \label{eq:T13}
      \displaystyle{q_{i}^{s} = \frac{1}{ \Delta x} \int_{x_{i-\frac{1}{2}}}^{x_{i+\frac{1}{2}}} q(x,\omega \Delta t) ) d x} \;.
\end{equation}
That is, the integrand $q(x,t) $ is evaluated at the fixed time $\omega \Delta t$, with $ 0\le \omega \le 1$.  To simplify the notation we have assumed that the time interval $[t_{n},t_{n+1}]$ has been translated to $[0,\Delta t]$.
For $\omega = \frac{1}{2}$ one obtains the second-order midpoint rule approximation.  For the linear source term (\ref{eq:RoeS3}), from  (\ref{eq:T13}) we obtain
\begin{equation}                           \label{eq:T14}
      q_{i}^{s} = \frac{1}{2} c  q_{i-1}^{n}+(1-\frac{1}{2} c) q_{i}^{n} \;.
\end{equation}
For the source term as in (\ref{eq:RoeS3}),  this is identical to (\ref{eq:T10b}),  with properties displayed in Table \ref{tab:BV}.

Another version,  applicable to the nonlinear case (\ref{eq:T2}),  uses the integral form of the conservation law to replace the right hand side of (\ref{eq:T13}) by
\begin{equation}                           \label{eq:T15}
     q_{i}^{s} =q_{i}^{n} -  \frac{\omega \Delta t}{\Delta x}( 
      f_{i+\frac{1}{2}} -  f_{i-\frac{1}{2}}) \;, 
\end{equation}
followed by (\ref{eq:T12}). That is,  the source term is $s_{i}  = s(q_{i}^{s})$. The complete scheme is
\begin{equation}                          \label{eq:T16}   
\left.
\begin{array}{ll}
       \mbox{Predictor:}  & q_{i}^{s} =q_{i}^{n} -  \frac{\omega \Delta t}{\Delta x}( 
      f_{i+\frac{1}{2}} -  f_{i-\frac{1}{2}}) \;,  \\[7pt]
       \mbox{Advection:}  & q_{i}^{adv} =q_{i}^{n} - \frac{\Delta t}{\Delta x}( 
      f_{i+\frac{1}{2}} -  f_{i-\frac{1}{2}}) \;, \\[7pt]
      \mbox{Corrector:}  & q_{i}^{n+1} = q_{i}^{adv} + \Delta t s(q_{i}^{s}) \;.
 \end{array}\right\}
\end{equation}
For the linear equation, scheme (\ref{eq:T16})  reproduces all the good upwind schemes proposed in \cite{Roe:1986c} and \cite{Bermudez:1994a}.   Moreover, in principle,  scheme (\ref{eq:T16}) extends to nonlinear inhomogeneous systems in a straightforward manner.   For the special cases $\omega=0$ (no upwinding of the source term) scheme (\ref{eq:T16}) reproduces the fully centred scheme, while for $\omega=1$ it reproduces the classical splitting scheme for source terms.

In the last 40 years, since Roe put forward his upwinding ideas to discretise source terms, the subject has undergone very significant advances.  The literature is vast and reviewing such advances lies outside the scope of these notes focused on Roe's early contributions.

\subsection{Nonlinear systems: Clarke's test}
\label{sec:NLsystems}

% DONE 5/9/26

Motivated by the numerical simulation of combustion
problems--characterised by the presence of reactive source
terms in the equations--John Clarke devised a test problem
with an exact solution \cite{Clarke:1987a},
\cite{Clarke:1985a}.  The time-dependent, one-dimensional
Euler equations for ideal gases were augmented by an
algebraic source term:
\begin{equation} \label{eq:Roe13}
  \left.
  \begin{array}{c}
    \partial_t {\bf Q}
    + \partial_x {\bf F}({\bf Q})
    = {\bf S}({\bf Q}) \;, \\[6pt]

    {\bf Q}
    =
    \left[
    \begin{array}{c}
      \rho \\
      \rho u \\
      E
    \end{array}
    \right] \;,

    {\bf F}({\bf Q})
    =
    \left[
    \begin{array}{c}
      \rho u \\
      \rho u^2 + p \\
      u(E+p)
    \end{array}
    \right] \;, \\[6pt]

    {\bf S}({\bf Q})
    \equiv H(-x)(C/a)
    \left[
    \begin{array}{c}
      \rho \\
      \rho u \\
      \rho h
    \end{array}
    \right] \;.
  \end{array}
  \right\}
\end{equation}
Here $\rho$ is density; $u$ is velocity; $p$ is pressure;
$E=\frac{1}{2}\rho u^{2}+\rho e$ is total energy;
$e=e(\rho,p)=p/[\rho(\gamma-1)]$ is specific internal
energy; $\gamma$ is the ratio of specific heats, taken here as 
$\gamma=1.4$;
$h=(E+p)/\rho$ is specific enthalpy; $H$ is the Heaviside
function; $C$ is a constant with dimensions of velocity;
and $a$ is the speed of sound.

The problem relates to the gas dynamics resulting from the
combustion of a high-energy reactive material after ignition,
as described by the 1D equations with appropriate source
terms (\ref{eq:Roe13}). The specific test problem considers
a cylindrical stick of reactive material inside a 1 m long
tube of constant cross-sectional area. At the initial time
$t=0$, ambient conditions are assumed throughout the domain.
An igniter vents hot gases of a specified chemical energy
uniformly on the left half of the cylinder. Once the
temperature has reached a prescribed value, the stick of
reactive material is assumed to ignite, thus generating a
gas-dynamical problem. For descriptions of the problem and
related physical and chemical parameters, see Clarke and
Toro \cite{Clarke:1985a}, Toro and Clarke
\cite{Toro:1985a}, Toro \cite{Toro:1986b}, and Clarke
\cite{Clarke:1987a}.

\begin{figure}[h]
      \centerline{
             \includegraphics[scale=0.10, angle=0]{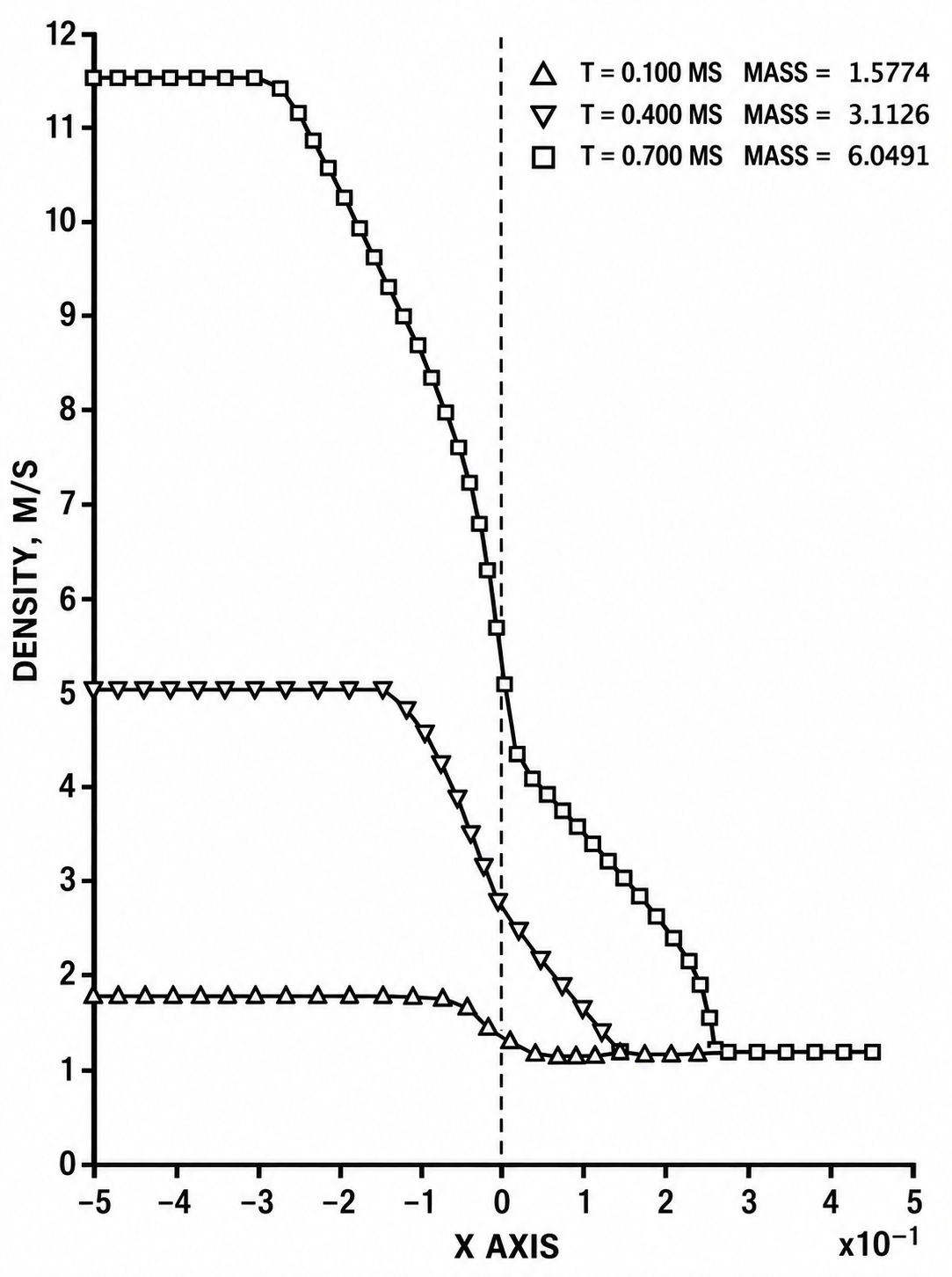}
            \includegraphics[scale=0.105, angle=0]{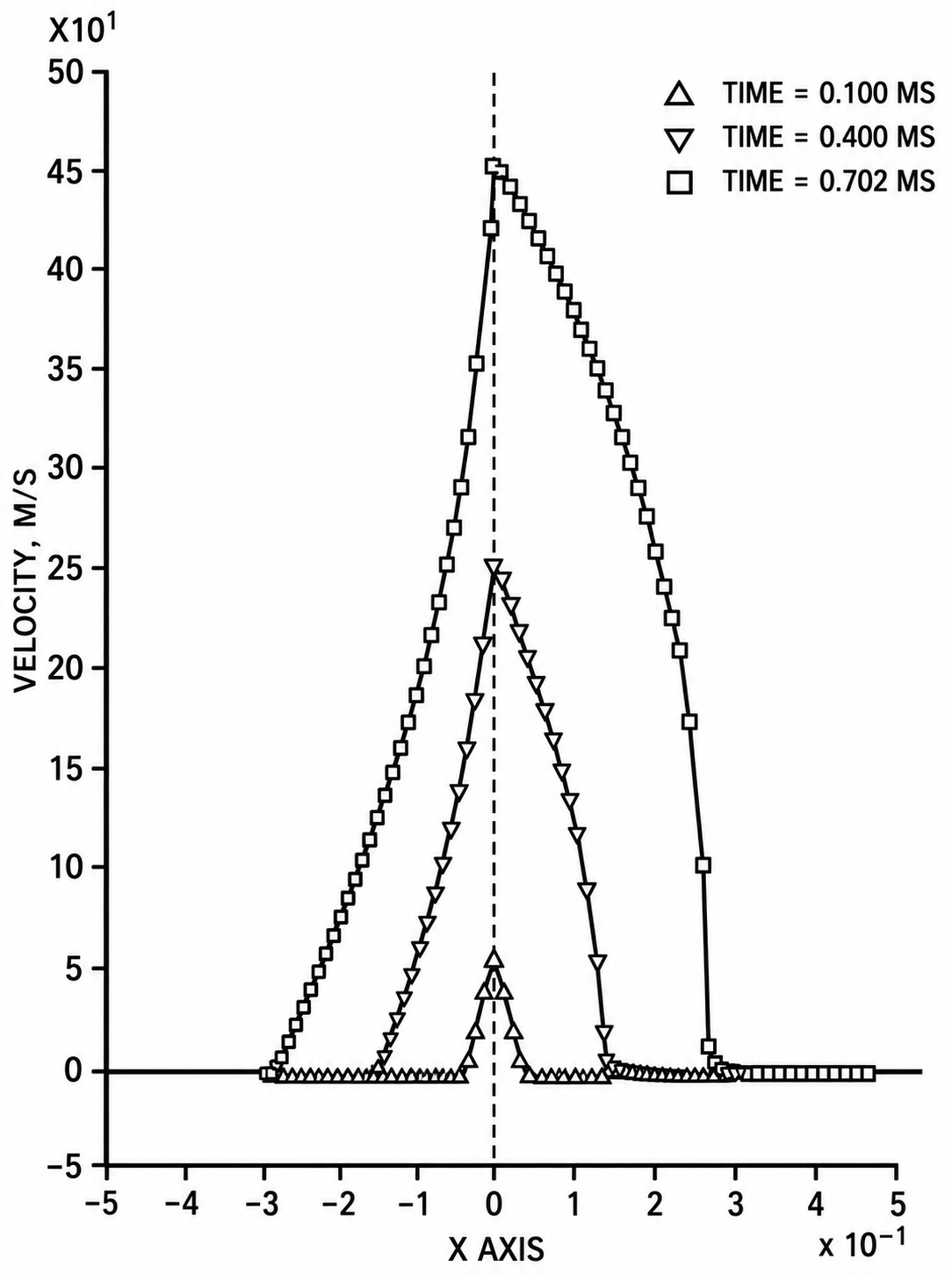}
      }
      \caption{Clarke's test problem.  Numerical solutions obtained from Roe's upwind 
      methodology for both the advective and source terms (symbols) and Clarke's exact 
      solution (line).  Results reproduced from \cite{Toro:1986b}.}
      \label{fig:ClarkeTest}
\end{figure}

The problem was simulated by solving equations
(\ref{eq:Roe13}) using a second-order flux-limited upwind
scheme with Roe's Riemann solver and upwinding for the source
term, as first proposed by Roe in \cite{Roe:1986c}.
Numerical results from the test problem and discussion are
found in \cite{Roe:1986c} and \cite{Toro:1986b}.
Fig.~\ref{fig:ClarkeTest} shows sample numerical results
using the same methodology, reproduced from
\cite{Toro:1986b}, with a mesh of 100 cells and the
SUPERBEE flux limiter. Numerical results (symbols) are
compared with Clarke's exact solution (line) at three output
times, for density (left frame) and velocity (right frame).
At a later time, the analytical solution becomes
multi-valued on the far right, so that the compression wave
on the right becomes a shock wave. Overall, the results are
quite satisfactory.

It is appropriate to mention that Roe had also kindly shared
his ideas with Paul Glaister, University of Reading, who
obtained satisfactory results for problems in cylindrical
and spherical symmetry using a similar numerical approach
\cite{Glaister:1985a}.

\section{R\&D at Cranfield}

%DONE 4/9/26

While still at Cranfield, Phil Roe contributed significantly to enhancing Cranfield's well-established relationship with industry through research, development and consultancy. Influenced by John Clarke's research interests in combustion waves and fluid dynamics, we became engaged in the modelling and simulation of detonation waves.

In the preliminary study of Clarke et al. \cite{Clarke:1989b}, we explored various computational strategies with the aim of identifying the best choice for computing problems involving strong detonation waves.

To this end, we adopted the detonation analogue devised by Fickett \cite{Fickett:1979a}, which exhibits many of the essential computational difficulties and possesses analytical solutions corresponding to overdriven and under-driven reacting flows.

We tested classical shock-capturing methods such as those of MacCormack and Godunov, as well as Roe's method \cite{Roe:1981a}, moving finite elements and Random Choice Methods (still in one dimension), especially using a new variant that yielded second-order accuracy \cite{Toro:1987a}.

In a second stage, in a project driven by PhD student James Quirk, the group developed realistic mathematical models for the theoretical simulation of unsteady, two-dimensional detonation waves in high-energy solids \cite{Quirk:1991a}, \cite{Clarke:1993a}. The model consisted of the two-dimensional compressible Euler equations augmented by extra conservation laws and source terms to account for chemical reactions and the tracking of material fronts. The numerical method selected for the full model was an extension of the conservative, shock-capturing technique of Roe \cite{Roe:1981a}.  In addition to the involvement of Clarke, Roe and Quirk, contributions to this project also came from  S. Karni,  L. G.  Simmonds and E. F.  Toro.

The treatment of reactive source terms in those years was challenging, as discussed at the time by Sweby \cite{Sweby:1989a} and Roe \cite{Roe:1986c}.

The tracking, or capturing, of material interfaces posed another numerical challenge. It was probably the first time that it was realised that systems in conservative form, including an additional advection equation to track the front, may produce oscillations in the vicinity of material interfaces. A most important additional ingredient of the work was an Adaptive Mesh Refinement (AMR) methodology that allowed the resolution of fine features such as reaction zones \cite{Clarke:1993a}.

Later works examined the potential of primitive and adaptive primitive/conservative methods \cite{Karni:1994a} to circumvent some of these difficulties. As a matter of fact, these failings of conservative methods are also present in multidimensional problems involving contact surfaces and shear waves \cite{Toro:2002b}.

\section{Roe's Way}

%DONE 2/8/26

In the abstract of his recent paper \cite{Roe:2020a}, Phil Roe states: {\it In this paper, the author recounts his forty-year plus struggle to find a sound basis for understanding the computational fluid dynamics of compressible flow.}

\begin{figure}[h]
      \centerline{
             \includegraphics[scale=0.3, angle=0]{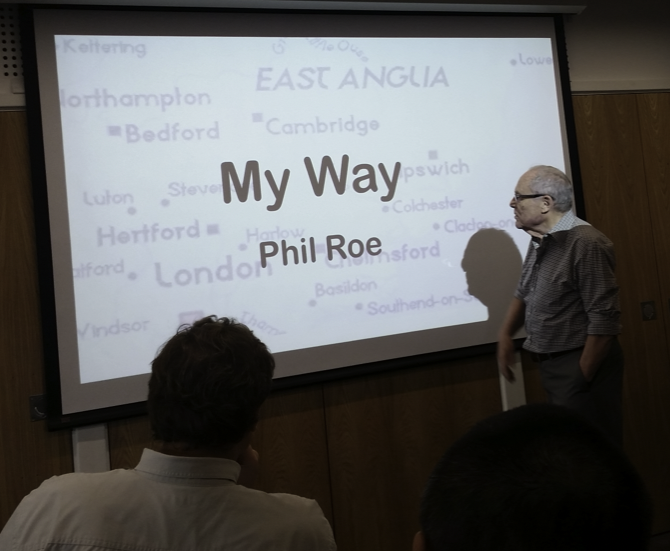}
      }
      \caption{Roe's 80th birthday celebration in 2018.  Phil  Roe  lecturing.  University of 
      Cambridge,  UK.}
      \label{fig:RoesWay}
\end{figure}

The statement and the subsequent content illustrate the long and unique road he travelled. The journey began in experimental aerodynamics,  driven by the very practical purpose of getting aircraft to fly.   It then followed a theoretical path for about five decades,  in search of better methods for solving the partial differential equations governing compressible fluid flow phenomena. The first phase--experimental research--strongly influenced the second phase, computational fluid dynamics. This influence was remarkably insightful, and perhaps unique,   owing to Roe's distinctive combination of experimental and theoretical experience.

The physical interpretation of every aspect of a computational method was always central to Roe's work.  And aircraft, of course.  Interestingly,  when recounting his move in 1984 from RAE to Cranfield, he states:

{\it This was another remarkable place: we had our own airfield...}

The airfield came first.

\subsection{Multidimensional upwinding}

% DONE

Upwinding was a key concept that drove Roe's research on numerical methods for partial differential equations associated with wave propagation phenomena. Beyond using the Riemann problem to devise upwind numerical fluxes, he also pioneered the adoption of upwinding in the discretisation of source terms \cite{Roe:1986c}. But perhaps the subject that spanned most of his research life was his desire to extend the concept of upwinding in a more complete manner than that offered by the one-dimensional Riemann problem.

Back in the Cranfield days in the 1980s, I vividly remember his triangle with arrows representing the directions of wave propagation. I confess I never took up the idea; it looked too challenging to me. But Phil was optimistic. I remember his talk at a workshop in Sendai (Japan) in 1991. He stated that, given three 
years, the subject of multidimensional upwinding would be complete.

At the same event, Ami Harten (1946--1994) was present. After hearing  Roe's list of the shortcomings of methods based on one-dimensional ideas, Harten stated that, really, the central problem was that of devising schemes with a high order of accuracy. 
Two towering figures in the field, with two different views.

But in fairness, Phil Roe was also sensitive to the issue of high accuracy, first in his work on limiters and TVD methods (high-order schemes), and later in the context of multidimensional upwinding, as seen in his later works.

When it came to technical issues,  Phil Roe had very strong views, and one could therefore be led to think that Phil was not appreciative of other people's ideas.  But this was far from the truth.  In fact, this is most clearly demonstrated by two of his most important lines of research concerned with multidimensional upwinding.

First, he appreciated the value of the early work by Ni in 1981 \cite{Ni:1981a}. In this work, Ni took the unknowns to be point values of the conserved variables and located them at the cell vertices; he then distributed the change among the four vertices. Ni's method became known as the Cell Vertex Finite Volume Method. Roe found that Ni's approach offered a suitable framework for his multidimensional upwinding doctrine.  Along with several of Roe's collaborators, the approach was further developed, giving rise to what became known as Residual Distribution Schemes 
\cite{Abgrall:2004a}, \cite{Ricchiuto:2010a}.

Another example is Roe's work on what is now known as Active Flux.  He started from van Leer's 1977 paper \cite{vanLeer:1977a} and chose scheme V.  This  scheme was virtually unknown in the community, but Phil discovered its merits. Thus came the task he set himself: to develop scheme V of van Leer \cite{Roe:2020a}.  On the distinguishing features of the method, literally,  Roe states: {\it First the interface fluxes are calculated from the initial data. This is done using the governing equations, rather than just by interpolation, so I have christened it the Active Flux scheme.}

Preliminary ideas on the Active Flux Method were communicated in 2011 by Eymann and Roe and later developed in Eymann's PhD thesis \cite{Eymann:2013a}. The Active Flux Method was applied to the Euler equations by Roe et al. \cite{Roe:2018a}. They compared results obtained with their scheme with those of discontinuous Galerkin finite element methods. In recent years, the Active Flux approach has drawn the attention of other research groups; see,  for example, Helzel et al. \cite{Helzel:2019a},  Barsukow \cite{Barsukow:2021a} and Abgrall and Liu \cite{Abgrall:2026a}.

Multidimensional upwinding was to remain Roe's passion for the rest of his life. In 2020 he wrote: {\it I spend much of my time these days preoccupied with multidimensional methods} \cite{Roe:2020a}.

\subsection{Phil Roe the philosopher}

%DONE

Phil Roe was, no doubt, a thinker, and therefore, in the truest sense, a philosopher; not a professional philosopher, of course. Conversing one-to-one with Phil was always a rewarding experience, in which his reflections and his sense of humour flowed smoothly and engagingly.  Often, at least for me,  the conversation would become intellectually challenging.  He was a good listener; you felt encouraged to develop your thoughts and opinions.  He often offered a synthesis of his own thoughts.

Here is a documented example of Phil's philosophical inclinations.  In discussing the interplay  between data (present) and the predicted solution (future) when designing a numerical method, Roe wrote:

{\it I am prepared to put forward a principle that can be rather generally applied: 

“Before attempting to predict the future, you must first understand the present”.}

Another instance is:

{\it I would like to propose the idea that the “order of difficulty” of a problem is the smallest number of good ideas that are needed to solve it. I think that we have learned by now that in CFD this number is greater than one.}

In his 2020 paper \cite{Roe:2020a},  Phil Roe wrote: {\it There is an aphorism, that I believe is due to the German philosopher Georg Christoph Lichtenberg, to the effect that ``knowledge does not mean all the things we happen to know, but only those we have thought about enough to know how they hang together and how they can be applied usefully''. I have a feeling that CFD has not yet reached that happy condition, and that the future may be just as exciting as the past.}

A good thinker is likely to hold strong convictions. Phil Roe was no exception: he held strong views on various issues close to his heart, particularly CFD and, more generally, numerical methods for partial differential equations. Engaging in robust discussions was not unusual for Phil. This trait contrasted somewhat with his mild and friendly social manner.

Probably around 1988--89, I accompanied him to a one-day meeting on CFD matters held at a UK government laboratory near Oxford. During this period, there was a strong drive in the UK to commercialise CFD methods for technological applications. Various individuals switched from academic research to applied, commercially relevant development and consultancy.  Phil had followed the opposite direction and had definite opinions about such people; in time 
some of them became known as ``CFD cowboys''.

Well, at that meeting near Oxford, having arrived a few minutes late, we sat right at the back of the lecture theatre, the highest point of the place. The organiser of the meeting was giving a one-hour lecture.  Soon I could feel the effect of the talk on Phil, sitting next to me and looking visibly uncomfortable in his seat. At question time, Phil asked the first question. The long answer  followed a winding path, so Phil interrupted the speaker to correct some of his statements. The second attempted answer was getting more and more complicated.  Phil Roe interrupted again to say {\it ....almost everything you say is wrong....} Even in the relative ignorance of a mere postdoctoral fellow, I felt Roe was right, technically, but I wouldn't have had the guts to be so forthright.

In 2006, at the Hyperbolic Conference in Lyon, France, Phil Roe engaged in a public debate with Sergei Godunov. These two giants had much in common, such as their emphasis on the physics of wave propagation in the design of numerical methods. But they also had strong differences of opinion.
I witnessed a repeat of such a debate in Saint Petersburg, Russia, in 2010. Public debate was not one of Godunov's strong points,  certainly not in English.

Phil Roe was also sensitive to more general issues of human concern.  One evening in a Cambridge college bar, I  witnessed an animated discussion, this time between two native English speakers. The subject was US environmental policies and politics; Roe's interlocutor was a rather strong-minded Chicago academic visiting Cambridge.

Empathy for the disadvantaged was part of his being; Phil was not merely an academic in his ivory tower, but a man profoundly conscious of the surrounding world and its miseries.

\section{Meeting Phil Again}

After Phil moved to Michigan in 1990 we remained in contact.   I myself moved away from Cranfield to Manchester in 1994 and then to Trento in Italy in 2002. We met at conferences on many occasions,  in many parts of the planet.  Places I remember offhand include  Michigan, Pennsylvania,  Sendai,  Manchester,  Bangalore, Santiago de Compostela, Oxford,  Beijing,  Saint Petersburg,  Lyon, Trento,  Cambridge,  M\'alaga and many more.  Often we would find an evening in which we could enjoy a pleasant conversation over dinner in some nearby restaurant.  And it was always a pleasure to talk with Phil.  We would talk about our families,  our children,  increasingly becoming independent.  And we could not avoid nostalgia and talk about the old Cranfield days.  Whenever he returned to Bedford for personal reasons we would occasionally manage to meet, socially.

Our chance meeting in Bangalore in 2005 was peculiar. I was taking part in the 25th International Symposium on Shock Waves,  held  at the Indian Institute of Science in Bangalore.  Having landed well after midnight, my senses that morning, during breakfast on the university campus, were not particularly alert. Slowly, slowly, however, while sipping my coffee with little enthusiasm, I began to pay attention to a conversation in English some distance away, behind me. I could not see the two men talking, but I recognised Phil's voice. Yes, it was unmistakably his voice, but I did not turn to look at the men. As they carried on talking, I got a pretty good idea of who the second man was. Then I turned. And there they were: Phil Roe and Bill Morton.
 
In 2011 Phil and his second wife Jacqueline were amongst the many friends and collaborators who joined my 65th birthday conference organised by my friend Elena V\'azquez in beautiful Santiago de Compostela,  Spain.

In 2016,  Phil and Jacqueline joined my 70th  birthday in Trento, Italy,   organised by my colleague and friend Michael Dumbser.

In 2018 my wife Brigitte and I joined Phil's 80th birthday conference in Cambridge,  organised by our common friend from the old Cranfield days, Nikos Nikiforakis.

My final interaction with Phil was remote and took place during a workshop at the Department of Mathematics at the Southern University of Science and Technology in Shenzhen, China, in December 2025. I was given the honour of introducing Phil’s presentation and moderating the discussion that followed. He delivered his talk online from Michigan. Once the connection was finally established, the presentation flowed very smoothly indeed, and Phil’s performance in the ensuing discussion was as sharp and engaging as ever.

This is the image of Phil Roe that I will preserve: still sharing his ideas with clarity and enthusiasm, and still engaging fully and generously in discussion.

\section{Acknowledgements}

For valuable feedback on the original manuscript,  I am grateful to Peter K. Sweby (University of Reading, UK), Smadar Karni (University of Michigan, USA),  Nikos Nikiforakis (University of Cambridge, UK), Elena V\'azquez (University of Santiago de Compostela, Spain),  Jacqueline Roe (University of Michigan, USA) and Alan Dawes (AWE, Aldermaston, UK).

\bibliographystyle{plain}
%\bibliography{/Users/toro/Documents/SharedReferences/refs-latest}
\bibliography{refs-latest}
\clearpage

%%%%%%%%%% PHOTOS
%
\begin{figure}
      \centerline{
             \includegraphics[scale=0.25, angle=0]{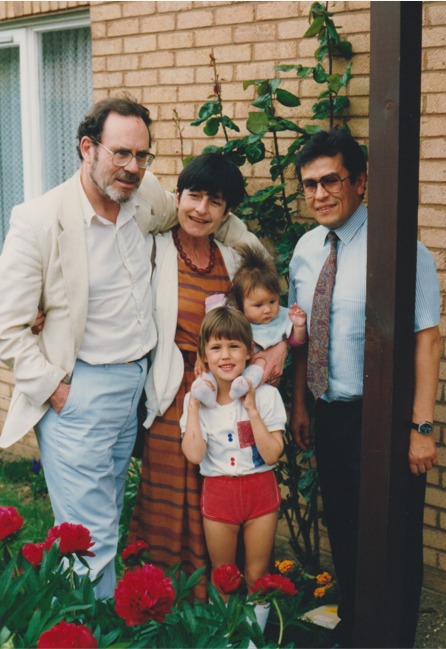}
      }
      \caption{Phil Roe (left), his first wife Susan,  Eleuterio Toro and his youngest 
      daughters Violeta and Eva,  at Toro's house in Bedford, 
      UK, 1992.}
      \label{fig:Roe-Kempston}
\end{figure}
\begin{figure}
      \centerline{
             \includegraphics[scale=0.15, angle=0]{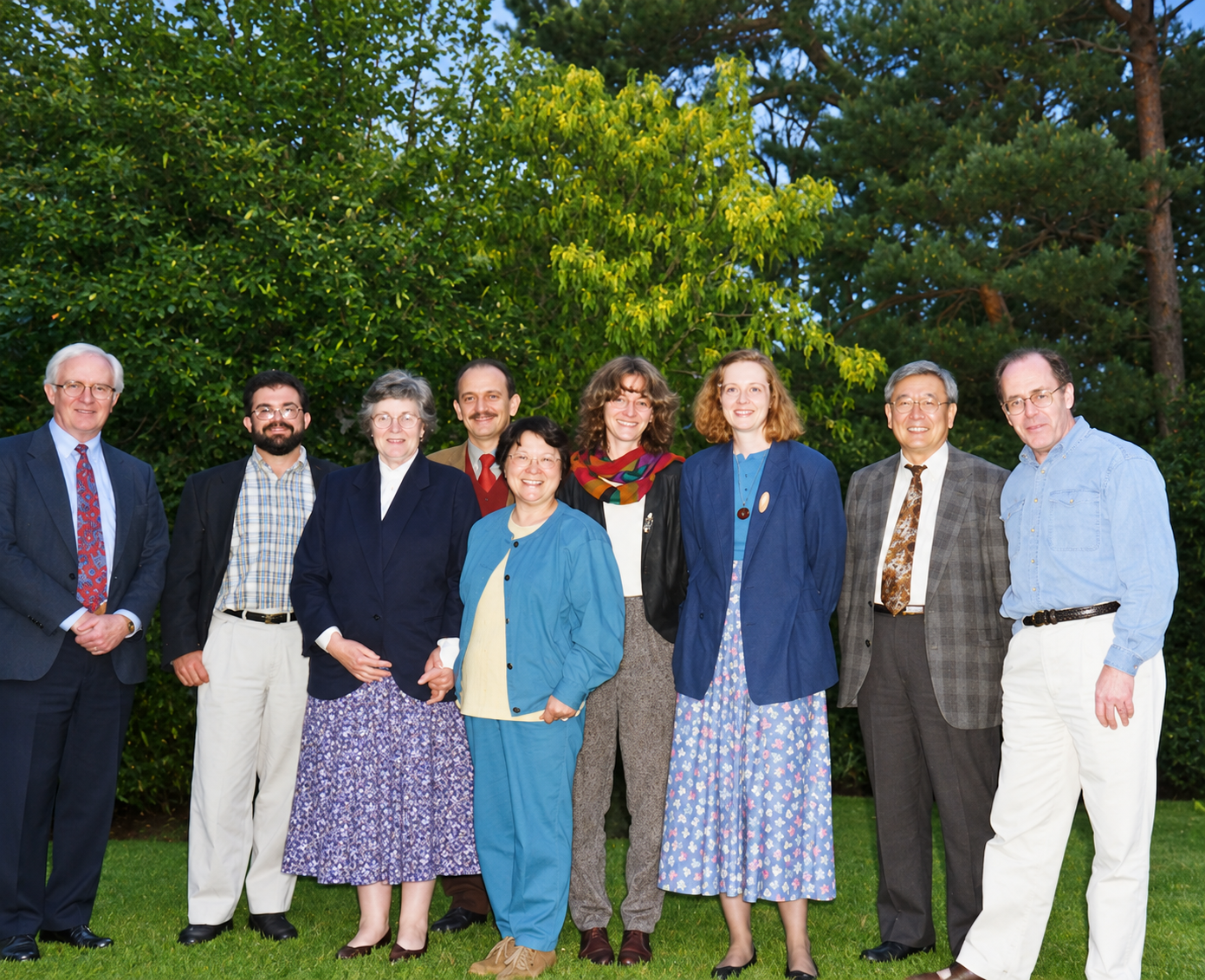}
      }
      \caption{Friends at my house in Manchester, May 1995, on occasion of the {\it International workshop on numerical methods for wave propagation phenomena. } From the right: Phil Roe, Kazuyoshi Takayama, Brigitte Toro, Smadar Karni, Chieko Takayama, Vincenzo Casulli, Jean Clarke, Jos\'e Mar\'ia  Mart\'i  and John Clarke.}
      \label{fig:Roe-Kempston}
\end{figure}
\begin{figure}
      \centerline{
             \includegraphics[scale=0.15, angle=0]{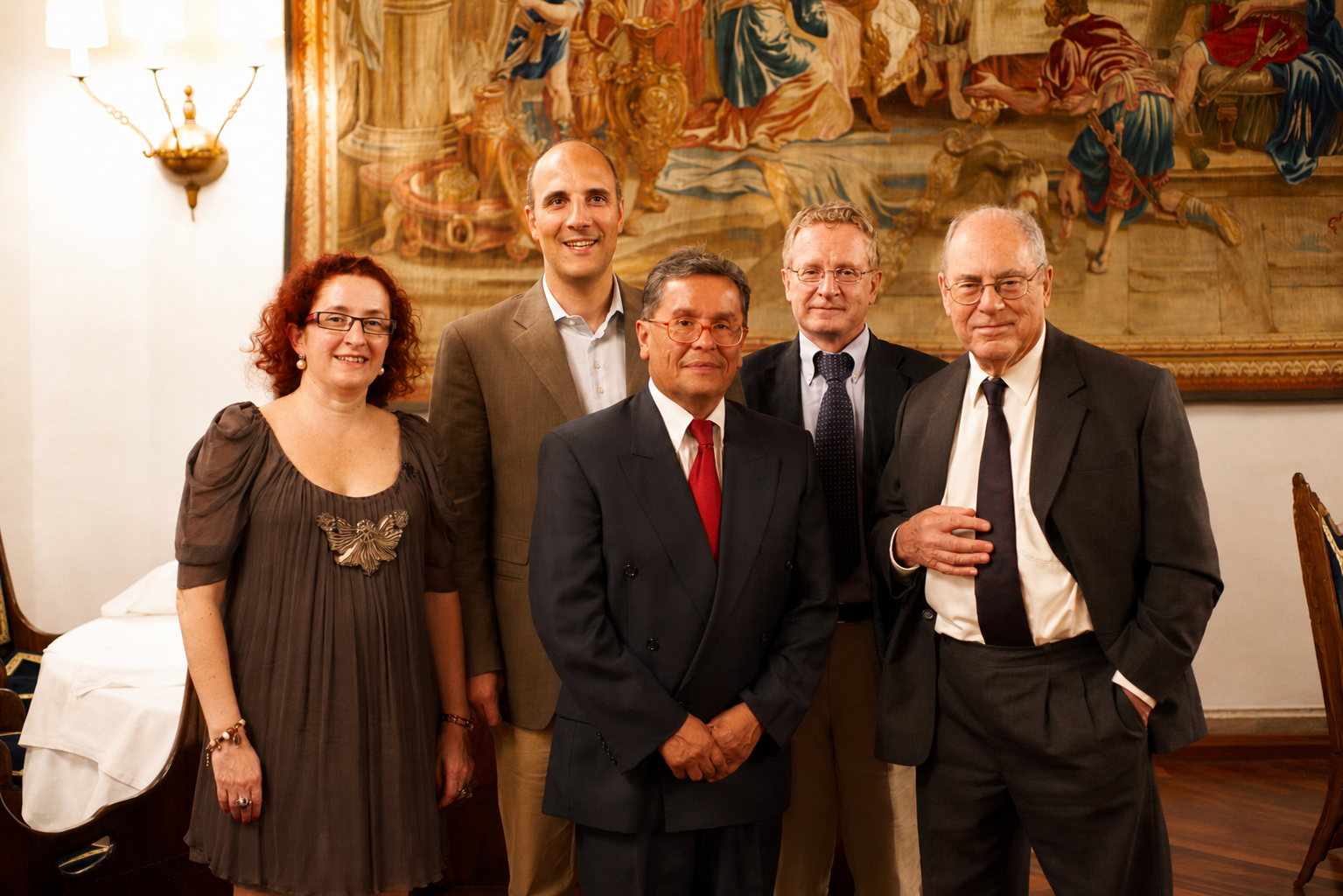}
      }
      \caption{Toro's 65th birthday celebration in 2011.  From the right: Phil Roe,  Aronne 
      Armanini, Eleuterio Toro,  Annunziato Siviglia and Elena V\'azquez.   Hostal de  
      los Reyes Cat\'olicos,  Santiago de Compostela, Spain.}
      \label{fig:PLRoe-USC}
\end{figure}
\begin{figure}
      \centerline{
             \includegraphics[scale=0.16, angle=0]{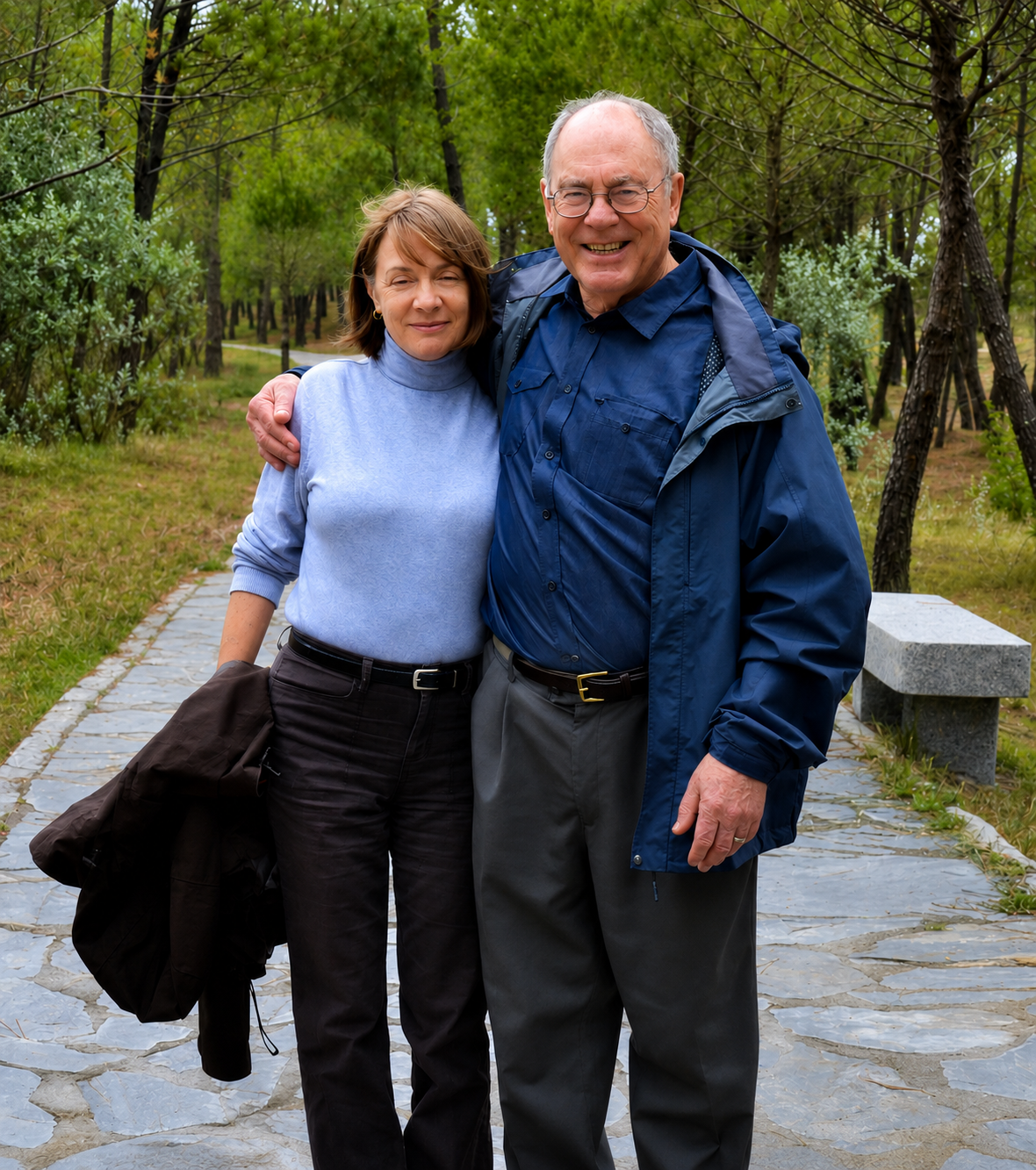}
      }
      \caption{Toro's 65th birthday celebration in 2011.  Phil Roe and his wife Jacqueline Roe.   
      Cape Finisterre,  Spain.}
      \label{fig:PLRoe-USC}
\end{figure}
\begin{figure}
      \centerline{
             \includegraphics[scale=0.15, angle=0]{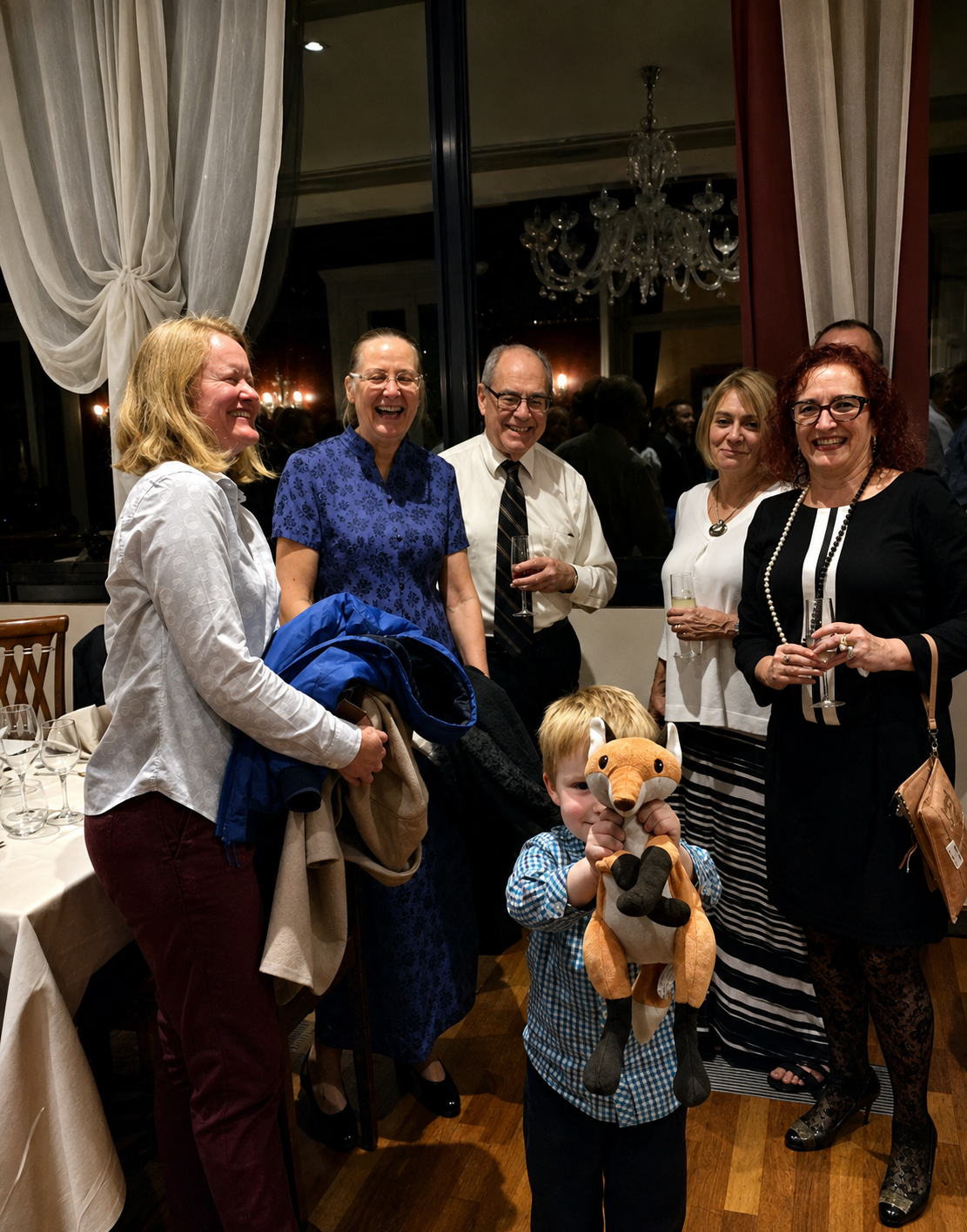}
      }
      \caption{Toro's 70th birthday celebration in 2016.  From the right: Elena 
      V\'azquez,  Jacqueline Roe,  Phil Roe,  Brigitte Toro, Lara Jamieson and Yiannis Nikiforakis. 
      University of Trento, Italy, 2016.}
      \label{fig:PLRoe-Trento}
\end{figure}
\begin{figure}
      \centerline{
             \includegraphics[scale=0.15, angle=0]{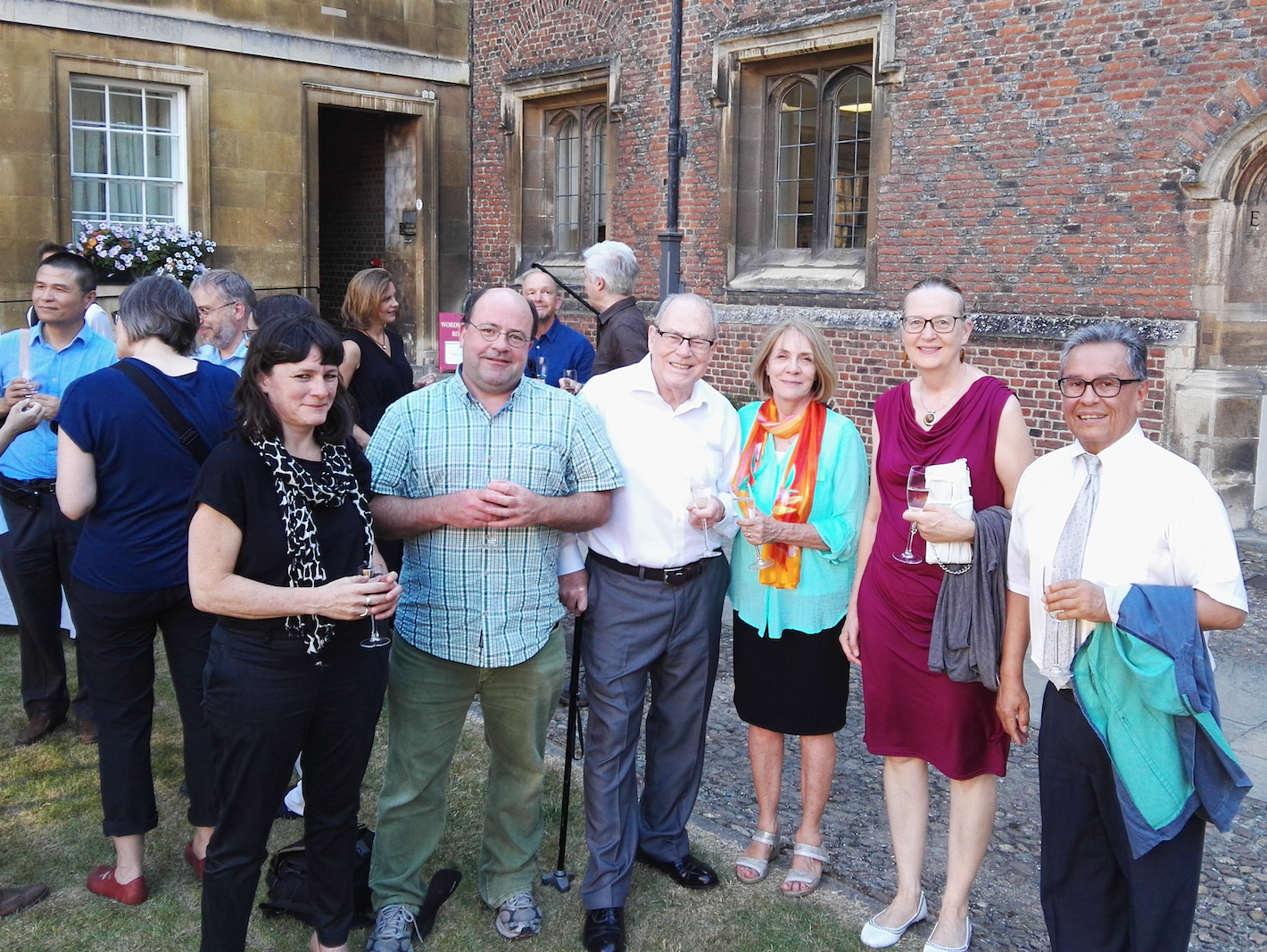}
      }
      \caption{Roe's 80th birthday celebration in 2018.  Phil Roe and his wife Jacqueline at the centre.  On the left Phil's  daughter Natasha and son Gerard.  On the right Eleuterio Toro and his wife Brigitte.  University of Cambridge,  UK.}
      \label{fig:PLRoe-Cranfield}
\end{figure}
\begin{figure}
      \centerline{
             \includegraphics[scale=0.13, angle=0]{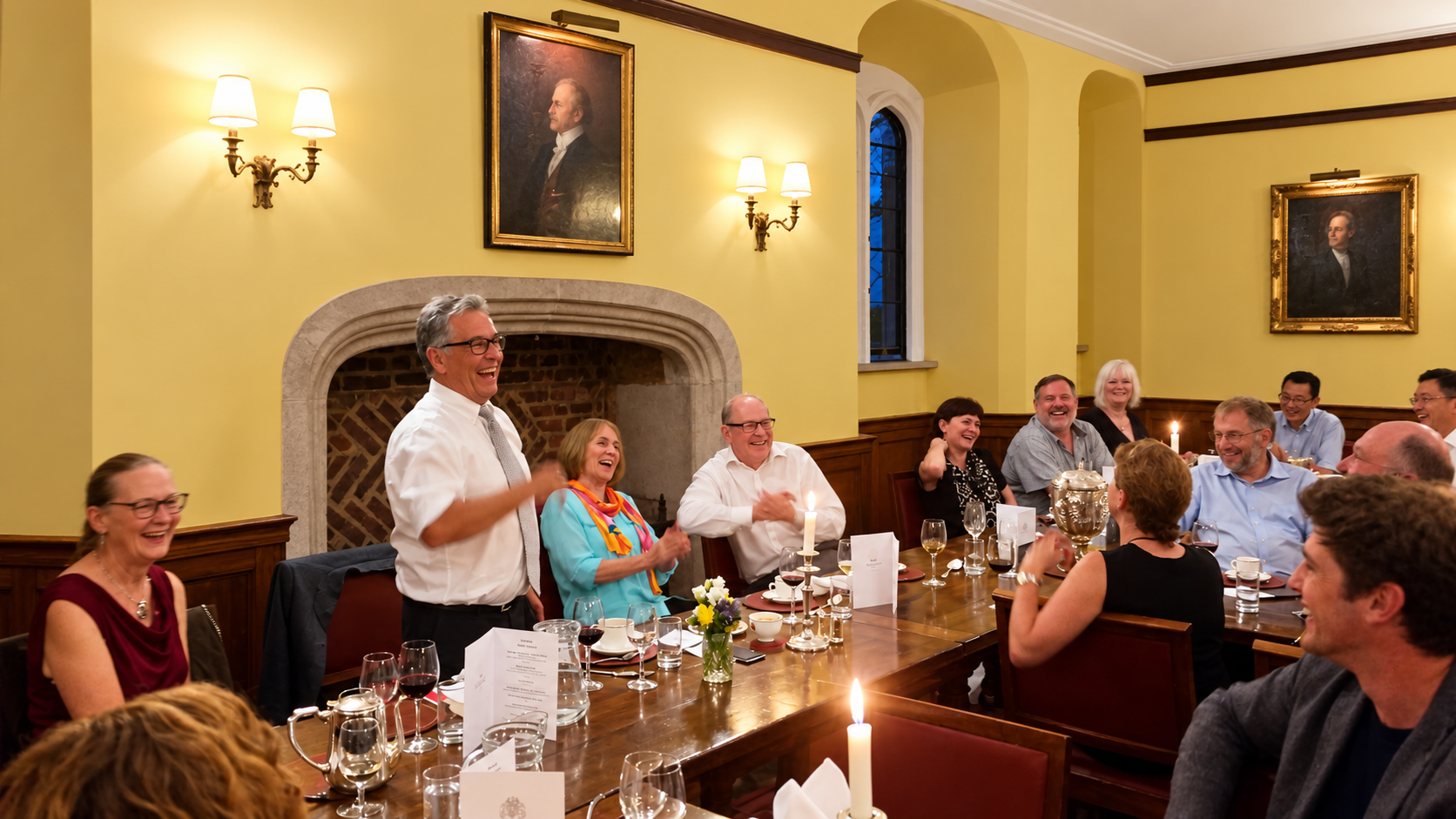}
      }
      \caption{ Roe's 80th birthday celebration in 2018.  Dinner at St John's College. From the left: Brigitte Toro, Eleuterio Toro,  Jacqueline Roe and Phil Roe.  University of Cambridge,  UK.}
      \label{fig:PLRoe-Cranfield}
\end{figure}

\end{document}